# Stable laws and domains of attraction in free probability theory

By Hari Bercovici and Vittorino Pata*
with an appendix by Philippe Biane

**Abstract**

In this paper we determine the distributional behavior of sums of free (in the sense of Voiculescu) identically distributed, infinitesimal random variables. The theory is shown to parallel the classical theory of independent random variables, though the limit laws are usually quite different. Our work subsumes all previously known instances of weak convergence of sums of free, identically distributed random variables. In particular, we determine the domains of attraction of stable distributions in the free theory. These freely stable distributions are studied in detail in the appendix, where their unimodality and duality properties are demonstrated.

## 1. Introduction

Denote by $\mathcal{M}$ the family of all Borel probability measures defined on the real line $\mathbf{R}$. Two measures $\mu, \nu$ in $\mathcal{M}$ will be said to be *equivalent* if there exist real numbers $a, b$, with $a > 0$, such that $\mu(S) = \nu(aS + b)$ for every Borel set $S \subset \mathbf{R}$. If $b = 0$ the equivalence relation will be said to be *strict*. We will write $\mu \sim \nu$ if $\mu$ and $\nu$ are equivalent. On the set $\mathcal{M}$ there are defined two associative composition laws denoted $*$ and $\boxplus$. The measure $\mu * \nu$ is the classical convolution of $\mu$ and $\nu$. In probabilistic terms, $\mu * \nu$ is the probability distribution of $X+Y$, where $X$ and $Y$ are (commuting) independent random variables with distributions $\mu$ and $\nu$, respectively. The measure $\mu \boxplus \nu$ is the free (additive) convolution of $\mu$ and $\nu$ introduced by Voiculescu [18] (for compactly supported measures; free convolution was extended by Maassen [11] to measures with finite variance and by Bercovici and Voiculescu [4] to the whole class $\mathcal{M}$). Thus, $\mu \boxplus \nu$ is the probability distribution of $X + Y$, where $X$ and $Y$ are free random variables with distributions $\mu$ and $\nu$, respectively. There are, naturally, free analogues of multiplicative convolution; these were first studied in [19].

---

*The first author was partially supported by a grant from the National Science Foundation.
1991 *Mathematics Subject Classification*. Primary 46L50, 60E07; Secondary 60E10.

An important class of measures occurs in connection with the study of the limit laws of probability. A measure $\nu \in \mathcal{M}$ will be said to be $*$-*infinitely divisible* if, for every natural number $n$, there exists a measure $\mu_n \in \mathcal{M}$ such that
$$\nu = \underbrace{\mu_n * \mu_n * \cdots * \mu_n}_{n \text{ times}}.$$
Analogously, a measure $\nu \in \mathcal{M}$ will be said to be $\boxplus$-*infinitely divisible* if, for every natural number $n$, there exists a measure $\mu_n \in \mathcal{M}$ such that
$$\nu = \underbrace{\mu_n \boxplus \mu_n \boxplus \cdots \boxplus \mu_n}_{n \text{ times}}.$$
Among the ($*$ or $\boxplus$)-infinitely divisible measures are the stable measures defined as follows. A measure $\nu$ which is not a point mass will be said to be $*$-*stable* (resp., $\boxplus$-*stable*) if for every $\mu_1, \mu_2 \in \mathcal{M}$ such that $\mu_1 \sim \nu \sim \mu_2$ it follows that $\mu_1 * \mu_2 \sim \nu$ (resp., $\mu_1 \boxplus \mu_2 \sim \nu$). Stability will be said to be *strict* if the equivalence relation $\sim$ is strict.

Measures which are $*$-infinitely divisible were studied by de Finetti, Kolmogorov, Lévy, and Hinčin. The $\boxplus$-infinitely divisible measures were introduced by Voiculescu [18] in the context of compact supports, and his results were extended in [11] and [4]. In particular, [4] contains a description of the $\boxplus$-stable laws, most of which do not have compact support.

Before outlining the role of infinitely divisible and stable measures in relation with limit laws we need a few more definitions. If $\nu_n$ and $\nu$ are elements of $\mathcal{M}$ or, more generally, finite Borel measures on $\mathbf{R}$, we say that $\nu_n$ converge to $\nu$ *weakly* if
$$\lim_{n \to \infty} \int_{-\infty}^{\infty} f(t)\, d\nu_n(t) = \int_{-\infty}^{\infty} f(t)\, d\nu(t)$$
for every bounded continuous function $f$ on $\mathbf{R}$. Fix now two measures $\mu, \nu \in \mathcal{M}$. We will say that $\mu$ belongs to the *partial $*$-domain of attraction* (resp., *partial $\boxplus$-domain of attraction*) of $\nu$ if there exist measures $\mu_1, \mu_2, \ldots$ equivalent to $\mu$, and natural numbers $k_1 < k_2 < \cdots$ such that
$$\underbrace{\mu_n * \mu_n * \cdots * \mu_n}_{k_n \text{ times}} \qquad \left(\text{resp., } \underbrace{\mu_n \boxplus \mu_n \boxplus \cdots \boxplus \mu_n}_{k_n \text{ times}}\right)$$
converges weakly to $\nu$ as $n \to \infty$. If the above convergence holds with $k_n = n$, we say that $\mu$ belongs to the $*$-*domain of attraction* (resp., $\boxplus$-*domain of attraction*) of $\nu$. We will denote by $\mathcal{P}_*(\nu)$ (resp., $\mathcal{D}_*(\nu)$, $\mathcal{P}_\boxplus(\nu)$, $\mathcal{D}_\boxplus(\nu)$) the partial $*$-domain of attraction (resp., $*$-domain of attraction, partial $\boxplus$-domain of attraction, $\boxplus$-domain of attraction) of $\nu$. The following result is due to Hinčin for the classical convolution. For free convolution it was proved in [13].

1.1. THEOREM. *A measure $\nu \in \mathcal{M}$ is $*$-infinitely divisible (resp., $\boxplus$-infinitely divisible) if and only if $\mathcal{P}_*(\nu)$ (resp., $\mathcal{P}_\boxplus(\nu)$) is not empty.*

The next theorem, perhaps the most important result of this work, gives an identification of the partial $\boxplus$-domains of attraction.

1.2. THEOREM. *There exists a bijection $\nu \leftrightarrow \nu'$ between $*$-infinitely divisible measures $\nu$ and $\boxplus$-infinitely divisible measures $\nu'$ such that $\mathcal{P}_*(\nu) = \mathcal{P}_\boxplus(\nu')$. More precisely, let $\mu_n \in \mathcal{M}$, let $k_1 < k_2 < \cdots$ be natural numbers, and set*

$$\nu_n = \underbrace{\mu_n * \mu_n * \cdots * \mu_n}_{k_n \text{ times}}, \qquad \nu'_n = \underbrace{\mu_n \boxplus \mu_n \boxplus \cdots \boxplus \mu_n}_{k_n \text{ times}}.$$

*Then $\nu_n$ converges weakly to $\nu$ if and only if $\nu'_n$ converges weakly to $\nu'$.*

This result shows that the large body of classical work dedicated to the study of limit laws for sums of independent identically distributed random variables can be transferred directly to the free context.

Passing now to stability, observe that a stable law belongs to its own domain of attraction, and in fact the following is true.

1.3. THEOREM. *Assume that $\nu \in \mathcal{M}$ is not a point mass. Then $\nu$ is $*$-stable (resp., $\boxplus$-stable) if and only if $\mathcal{D}_*(\nu)$ (resp., $\mathcal{D}_\boxplus(\nu)$) is not empty.*

The definition of stable laws, and the above result, could be extended to include point masses, but there are good reasons for not doing so. Convergence to a point mass is the object of laws of large numbers.

Theorem 1.3 was proved by Lévy [9] in the classical case and by Pata [12] in the free case.

The main result we prove for stability is as follows. This result was announced in [1], but the proof we present here is not the one envisaged at the time [1] was written.

1.4. THEOREM. *A measure $\mu \in \mathcal{M}$ belongs to a $*$-domain of attraction if and only if it belongs to a $\boxplus$-domain of attraction. More precisely, there exists a bijection $\nu \leftrightarrow \nu'$ between $*$-stable laws $\nu$ and $\boxplus$-stable laws $\nu'$ such that $\mathcal{D}_*(\nu) = \mathcal{D}_\boxplus(\nu')$.*

A particular case of this result, related to the central limit theorem, was proved earlier. Namely, if $X_1, X_2, \ldots$ are nonconstant, bounded, free identically distributed random variables with mean zero and variance one, it was shown by Voiculescu [17] that

$$\frac{X_1 + X_2 + \cdots + X_n}{\sqrt{n}} \xrightarrow[n \to \infty]{} \omega$$

in distribution, where $\omega$ is the semicircle law with density

$$d\omega(t) = \begin{cases} \dfrac{1}{2\pi}\sqrt{4-t^2}\,dt & \text{if } t \in [-2, 2], \\ 0 & \text{otherwise.} \end{cases}$$

The semicircle law $\omega$ first appeared prominently as a weak limit in Wigner's work on the eigenvalues of large random matrices [22]. The free central limit theorem was extended to unbounded variables with finite variance by Maassen [11]. Finally, Pata [14] proved the following result, which we can now view as a particular case of Theorem 1.4.

1.5. THEOREM. *Denote by $\nu$ the standard Gaussian (or normal) distribution and by $\nu'$ the standard semicircle distribution. Then we have $\mathcal{D}_*(\nu) = \mathcal{D}_\boxplus(\nu')$.*

An analogous result for the weak law of large numbers was established in [10] and [2]. Namely, it was shown that the classical weak law of large numbers holds for a measure $\mu$ if and only if the free weak law of large numbers holds for $\mu$. This result can now be seen as an immediate consequence of Theorem 1.2 because the Dirac measure at the origin $\nu = \delta_0$ corresponds with $\nu' = \delta_0$. The necessary and sufficient condition for the classical weak law of large numbers was given by Kolmogorov in terms of tail-sums. It is as follows:

$$\lim_{t\to\infty} t\mu(\{x : |x| > t\}) = 0.$$

The correspondence $\nu \leftrightarrow \nu'$ is the same in Theorems 1.2 and 1.4. This correspondence is easily seen from the classical Lévy-Hinčin formula and its free counterpart (see [18], [11], and [4]). This may seem to be a rather formal correspondence, but Theorem 1.2 shows the relationship to be quite deep.

There is a third operation on $\mathcal{M}$ denoted $\uplus$ and called the Boolean convolution. This was introduced by Speicher and Woroudi [16]. We will show in Section 6 that our results also extend to Boolean convolution. Infinitely divisible and (strictly) stable measures were determined in [16], where some of the basic limit theorems were proved, e.g., the central limit theorem for measures with finite variance. Curiously, all measures in $\mathcal{M}$ are infinitely divisible relative to $\uplus$.

The appendix contributed by Philippe Biane gives a description of the densities of $\boxplus$-stable distributions. This description allows him to prove that these distributions are unimodal and satisfy a duality relation, analogous to Zolotarev's relation satisfied by the $*$-stable distributions.

The remainder of this paper is organized as follows. In Section 2 we discuss the calculation of free convolutions and we develop an asymptotic formula for Voiculescu's analogue of the Fourier transform which is the key to our results. Section 3 is dedicated to the proof of Theorem 1.2. The results about domains of attraction are deduced in Section 4. In Section 5 we give a description of the Cauchy transforms of measures belonging to a domain of attraction. The corresponding limit theorems can also be based on the formulas given here, and this is the approach we described briefly in [1].

## 2. Cauchy transforms and free convolution

As mentioned in the introduction, the free convolution $\mu \boxplus \nu$ of two measures $\mu, \nu \in \mathcal{M}$ is the distribution of $X + Y$, where $X$ and $Y$ are free random variables with distributions $\mu$ and $\nu$, respectively. We will not enter the details of free random variables, but refer the reader to the monograph [21], and to [4] for the discussion of unbounded random variables. We concentrate on the analytic apparatus needed for the calculation of free convolutions.

Denote by $\mathbf{C}$ the complex plane and set $\mathbf{C}^+ = \{z \in \mathbf{C} : \Im z > 0\}$, $\mathbf{C}^- = -\mathbf{C}^+$. For a measure $\mu \in \mathcal{M}$ one defines the *Cauchy transform* $G_\mu : \mathbf{C}^+ \to \mathbf{C}^-$ by

$$G_\mu(z) = \int_{-\infty}^{\infty} \frac{1}{z - t} \, d\mu(t), \quad z \in \mathbf{C}^+.$$

We also set $F_\mu(z) = 1/G_\mu(z)$, $z \in \mathbf{C}^+$, so that $F_\mu : \mathbf{C}^+ \to \mathbf{C}^+$ is analytic. It is easy to see that $F_\mu(z)/z$ tends to 1 as $z \to \infty$ nontangentially to $\mathbf{R}$ (i.e., such that $\Re z/\Im z$ stays bounded), and this implies that $F_\mu$ has certain invertibility properties. To be precise, for two numbers $\eta, M > 0$ we set

$$\Gamma_\eta = \{z = x + iy \in \mathbf{C}^+ : |x| < \eta y\}$$

and

$$\Gamma_{\eta,M} = \{z = x + iy \in \Gamma_\eta : y > M\}.$$

Then for every $\eta > 0$ there exists $M = M(\mu, \eta)$ such that $F_\mu$ has a left inverse $F_\mu^{-1}$ defined on $\Gamma_{\eta,M}$. The function

$$\phi_\mu(z) = F_\mu^{-1}(z) - z$$

will be called the *Voiculescu transform* of $\mu$. It is not hard to show that $\Im \phi_\mu(z) \leq 0$ for $z$ in a truncated cone $\Gamma_{\eta,M}$ where $\phi_\mu$ is defined. We also have $\lim \phi_\mu(z)/z = 0$ as $|z| \to \infty$, $z \in \Gamma_\eta$; in other words, $\phi_\mu(z) = o(z)$ as $|z| \to \infty$, $z \in \Gamma_\eta$.

2.1. THEOREM. *If $\mu, \nu \in \mathcal{M}$, then $\phi_{\mu \boxplus \nu} = \phi_\mu + \phi_\nu$ in any truncated cone $\Gamma_{\eta,M}$ where all three functions involved are defined.*

This remarkable result was discovered by Voiculescu [18]. He only considered compactly supported measures $\mu$, in which case $\phi_\mu$ is defined in a neighborhood of $\infty$. The result was extended by Maassen [11] to measures with finite variance; the general case was proved in [4].

2.2. LEMMA. *Let $\mu \in \mathcal{M}$, $n$ a natural number, and set*

$$\mu_n = \underbrace{\mu \boxplus \mu \boxplus \cdots \boxplus \mu}_{n \text{ times}}.$$

*Assume that $F_{\mu_n}^{-1}$ exists in a truncated cone $\Gamma_{\eta,M}$. Then $F_\mu^{-1}$ also exists in that truncated cone.*

*Proof.* Let $\eta' < \eta$ and $M' > M$ be such that $F_\mu^{-1}$ is defined in $\Gamma_{\eta',M'}$ and the relationship $\phi_{\mu_n} = n\phi_\mu$ holds in $\Gamma_{\eta',M'}$. An easy calculation shows that

$$F_\mu^{-1}(z) = \frac{1}{n} F_{\mu_n}^{-1}(z) + \left(1 - \frac{1}{n}\right) z, \quad z \in \Gamma_{\eta',M'}.$$

Now, the right-hand side $\psi(z) = F_{\mu_n}^{-1}(z)/n + (1-1/n)\,z$ is defined in $\Gamma_{\eta,M}$, it takes values in $\mathbf{C}^+$, and

$$F_\mu(\psi(z)) = z, \quad z \in \Gamma_{\eta',M'}.$$

By the uniqueness of analytic continuation we must have $F_\mu(\psi(z)) = z$ for all $z \in \Gamma_{\eta,M}$; hence $F_\mu^{-1}$ is defined in $\Gamma_{\eta,M}$ and is equal to $\psi$ there. $\square$

In the study of limit laws it is important to translate weak convergence of probability measures into convergence properties of the corresponding Voiculescu transforms. This is achieved by the following result (cf. [4] and [2]).

2.3. PROPOSITION. *Let $\mu_n \in \mathcal{M}$ be a sequence. The following assertions are equivalent*:
(i) *$\mu_n$ converges weakly to a probability measure $\mu$;*
(ii) *there exist $\eta, M > 0$ such that the sequence $\phi_{\mu_n}$ converges uniformly on $\Gamma_{\eta,M}$ to a function $\phi$, and $\phi_{\mu_n}(z) = o(|z|)$ uniformly in $n$ as $z \to \infty$, $z \in \Gamma_{\eta,M}$;*
(iii) *there exist $\eta', M' > 0$ such that the functions $\phi_{\mu_n}$ are defined on $\Gamma_{\eta',M'}$ for every $n$, $\lim_{n\to\infty} \phi_{\mu_n}(iy)$ exists for every $y > M'$, and $\phi_{\mu_n}(iy) = o(y)$ uniformly in $n$ as $y \to \infty$.*

*Moreover, if (i) and (ii) are satisfied, then $\phi = \phi_\mu$ in $\Gamma_{\eta,M}$.*

We will require a result about analytic functions which was already noted in a different form in [14].

2.4. LEMMA. *Let $\eta, M, \varepsilon$ be positive numbers, and let $\phi : \Gamma_{\eta,M} \to \mathbf{C}$ be an analytic function such that*

$$|\phi(z)| \le \varepsilon|z|, \quad z \in \Gamma_{\eta,M}.$$

*For every $\eta' < \eta$ and $M' > M$ there exists $k > 0$ such that*

$$|\phi'(z)| \leq k\varepsilon, \quad z \in \Gamma_{\eta',M'}.$$

*Proof.* Fix $z \in \Gamma_{\eta,M}$, and let $r = r(z)$ be the largest number such that the circle $\{\zeta : |\zeta - z| = r\}$ is contained in the closure of $\Gamma_{\eta,M}$. Cauchy's inequality then yields

$$|\phi'(z)| \leq \frac{1}{r} \sup_{|\zeta-z|=r} \varepsilon|\zeta| = \frac{\varepsilon}{r}(|z| + r) \leq \frac{2\varepsilon|z|}{r}.$$

Now, $r$ can be calculated explicitly:

$$r = \min\left\{\frac{\eta y - |x|}{\sqrt{1+\eta^2}}, y - M\right\}$$

if $z = x + iy$. It is clear that

$$k = \sup_{z \in \Gamma_{\eta',M'}} \frac{2|z|}{r(z)}$$

is finite if $\eta' < \eta$ and $M' > M$. □

In what follows, the statement "$z \to \infty$ nontangentially" will simply mean that $|z| \to \infty$ but $\Re z / \Im z$ stays bounded, i.e., $z$ stays within a cone $\Gamma_\eta$. Lemma 2.4 implies that if $\phi(z) = o(z)$ as $z \to \infty$ nontangentially, then $\phi'(z) = o(1)$ as $z \to \infty$ nontangentially.

The actual calculation of $\phi_\mu$ requires inverting an analytic function. In order to avoid this inversion we will develop an approximation for $\phi_\mu$ which will allow us to formulate weak limit theorems directly in terms of Cauchy transforms.

2.5. PROPOSITION. *For every $\mu \in \mathcal{M}$ we have*

$$\phi_\mu(z) = z^2 \left[G_\mu(z) - \frac{1}{z}\right](1 + o(1))$$

*as $z \to \infty$ nontangentially.*

*Proof.* We have $F_\mu(z) = z + u_\mu(z)$ with $u_\mu(z)/z \to 0$ as $z \to \infty$ nontangentially. An application of Lemma 2.4 (and the remark following it) implies that $u'_\mu(z) = o(1)$ as $z \to \infty$ nontangentially. Thus we have $F'_\mu(z) = 1 + o(1)$ as $z \to \infty$ nontangentially, and this also implies that $(F_\mu^{-1})'(z) = 1 + o(1)$ as $z \to \infty$ nontangentially. If $|z|$ is sufficiently large then the straight line segment joining $z$ and $F_\mu(z)$ lies in a truncated cone where $\phi_\mu$ is defined. Thus

$$\begin{aligned}
\phi_\mu(z) &= F_\mu^{-1}(z) - z \\
&= F_\mu^{-1}(z) - F_\mu^{-1}(F_\mu(z)) \\
&= \int_{F_\mu(z)}^z (F_\mu^{-1})'(\zeta)\, d\zeta \\
&= \int_{F_\mu(z)}^z (1 + o(1))\, d\zeta \\
&= (z - F_\mu(z))(1 + o(1))
\end{aligned}$$

as $z \to \infty$ nontangentially. To conclude the proof observe that

$$z - F_\mu(z) = z^2 \left[ G_\mu(z) - \frac{1}{z} \right] \frac{1}{zG_\mu(z)},$$

and $zG_\mu(z) = 1 + o(1)$ as $z \to \infty$ nontangentially. $\square$

The preceding statement can be made uniform over a tight family of measures $\mu$. Recall that a family $\mathcal{F} \subset \mathcal{M}$ is *tight* if

$$\lim_{N \to \infty} \sup_{\mu \in \mathcal{F}} \mu(\{t : |t| > N\}) = 0.$$

2.6. PROPOSITION. *Let $\mathcal{F} \subset \mathcal{M}$ be a tight family.*

(1) *For every $\eta > 0$ there exists $M > 0$ such that $F_\mu^{-1}$ (and hence $\phi_\mu$) is defined in $\Gamma_{\eta, M}$ for every $\mu \in \mathcal{F}$.*
(2) *Let $\eta$ and $M$ be such that $F_\mu^{-1}$ is defined in $\Gamma_{\eta, M}$ for every $\mu \in \mathcal{F}$, and write*

$$\phi_\mu(z) = z^2 \left[ G_\mu(z) - \frac{1}{z} \right] (1 + v_\mu(z)), \quad z \in \Gamma_{\eta, M}, \mu \in \mathcal{F}.$$

*Then $\sup_{\mu \in \mathcal{F}} |v_\mu(z)| = o(1)$ as $z \to \infty$ nontangentially.*

*Proof.* Define the functions $u_\mu$ as in the proof of Proposition 2.5. From that proof (and from the uniform nature of Lemma 2.4) we see that we only need to prove that $u_\mu(z)/z \to 0$ uniformly in $\mu \in \mathcal{F}$ as $z \to \infty$ nontangentially. To do this, fix $N > 0$ and $z = x + iy \in \mathbf{C}^+$. We have

$$\left| \frac{t}{z - t} \right| \leq \sqrt{1 + \left( \frac{x}{y} \right)^2}, \quad t \in \mathbf{R},$$

with equality for $t = (x^2 + y^2)/x$ if $x \neq 0$, and

$$\left| \frac{t}{z - t} \right| \leq \frac{N}{y}, \quad t \in (-N, N).$$

Thus for every $\mu \in \mathcal{M}$,

$$|zG_\mu(z) - 1| = \left|\int_{-\infty}^\infty \left(\frac{z}{z-t} - 1\right) d\mu(t)\right|$$

$$\leq \int_{-\infty}^\infty \left|\frac{t}{z-t}\right| d\mu(t)$$

$$\leq \frac{N}{y} + [1 - \mu((-N, N))]\sqrt{1 + \left(\frac{x}{y}\right)^2}.$$

When the last quantity above is less than $1/2$ we can also make the following estimate:

$$\left|\frac{F_\mu(z)}{z} - 1\right| = \left|\frac{zG_\mu(z) - 1}{zG_\mu(z)}\right|$$

$$\leq \frac{|zG_\mu(z) - 1|}{1 - |zG_\mu(z) - 1|}$$

$$\leq 2\left\{\frac{N}{y} + [1 - \mu((-N, N))]\sqrt{1 + \left(\frac{x}{y}\right)^2}\right\}.$$

Assume now that we are given $\varepsilon, \eta > 0$, and $z \in \Gamma_\eta$. Then $N$ can be chosen so that

$$[1 - \mu((-N, N))]\sqrt{1 + \eta^2} < \frac{\varepsilon}{4}$$

for every $\mu \in \mathcal{F}$. We can also choose $M$ so that $N/M < \varepsilon/4$, and we deduce that for $z = x + iy \in \Gamma_{\eta,M}$,

$$\left|\frac{F_\mu(z)}{z} - 1\right| \leq \varepsilon$$

for every $\mu \in \mathcal{F}$ (provided that $\varepsilon < 1/2$). This demonstrates the uniformity of the convergence of $u_\mu$ to zero. $\square$

A useful related result is as follows.

2.7. PROPOSITION. *Let $\mu_n \in \mathcal{M}$ be a sequence converging weakly to $\delta_0$. Then there exist $\eta > 0$ and $M > 0$ such that $\mu_n$ has a Voiculescu transform of the form*

$$\phi_{\mu_n}(z) = z^2\left[G_{\mu_n}(z) - \frac{1}{z}\right](1 + v_n(z)), \quad z \in \Gamma_{\eta,M},$$

*and*

$$\lim_{n\to\infty} v_n(z) = 0, \quad z \in \Gamma_{\eta,M}.$$

*Proof.* As in Proposition 2.6, it is possible to find a truncated cone $\Gamma_{\eta,M}$ such that $\phi_{\mu_n}(z)$ has the above representation for $z \in \Gamma_{\eta,M}$. We must only prove the last assertion. Upon restricting the truncated cone, we know from Proposition 2.3 that $\lim_{n\to\infty} \phi_{\mu_n} = 0$ uniformly on the bounded subsets of $\Gamma_{\eta,M}$, and $\phi_{\mu_n}(z) = o(|z|)$ uniformly in $n$ as $z \to \infty$, $z \in \Gamma_{\eta,M}$. Fix now $z \in \Gamma_{\eta,M}$. Then, in particular, $z \in \Gamma_{\eta',M'}$ for some $\eta' < \eta$ and $M' > M$. Select $\varepsilon > 0$. It is clear that, for every $\omega \in \Gamma_{\eta',M'}$, we have $|\phi_{\mu_n}(\omega)| < \varepsilon|\omega|$, for $n$ big enough. Thus Lemma 2.4 entails $|\phi'_{\mu_n}(w)| < k\varepsilon$, for some fixed $k > 0$ (independent of $\varepsilon$) for $n$ big enough. Since $\mu_n$ converges weakly to $\delta_0$, it is also clear that $F_{\mu_n}(z) \underset{n\to\infty}{\longrightarrow} z$, and therefore, for $n$ big enough, the segment joining $z$ and $F_{\mu_n}(z)$ is entirely contained in $\Gamma_{\eta',M'}$. Proceeding as in the proof of Proposition 2.5, we have

$$\phi_{\mu_n}(z) = z - F_{\mu_n}(z) + \int_{F_{\mu_n}(z)}^{z} \phi'_{\mu_n}(\zeta)\, d\zeta$$

$$= z^2 \left[G_{\mu_n}(z) - \frac{1}{z}\right] + z^2 \left[G_{\mu_n}(z) - \frac{1}{z}\right] \left[\frac{1}{zG_{\mu_n}(z)} - 1\right]$$

$$+ \int_{F_{\mu_n}(z)}^{z} \phi'_{\mu_n}(\zeta)\, d\zeta.$$

The proof is completed by observing that

$$\left|\int_{F_{\mu_n}(z)}^{z} \phi'_{\mu_n}(\zeta)\, d\zeta\right| \leq |z - F_{\mu_n}(z)|k\varepsilon,$$

and $\lim_{n\to\infty} 1/zG_{\mu_n}(z) = 1$. $\square$

### 3. The Lévy-Hinčin formula and limit laws

A characterization of $*$-infinitely divisible measures $\nu \in \mathcal{M}$ is given by the well-known Lévy-Hinčin formula in terms of the Fourier transform (or characteristic function)

$$(\mathcal{F}\nu)(t) = \int_{-\infty}^{\infty} e^{itx}\, d\nu(x).$$

3.1. THEOREM. *A measure $\nu \in \mathcal{M}$ is $*$-infinitely divisible if and only if there exist a finite positive Borel measure $\sigma$ on $\mathbf{R}$, and a real number $\gamma$ such that*

$$(\mathcal{F}\nu)(t) = \exp\left[i\gamma t + \int_{-\infty}^{\infty} (e^{itx} - 1 - itx)\frac{x^2+1}{x^2}\, d\sigma(x)\right], \quad t \in \mathbf{R},$$

*where $(e^{itx} - 1 - itx)(x^2+1)/x^2$ must be interpreted as $-t^2/2$ for $x = 0$.*

The free situation was settled in the general case in [4] (see also [3]), where it was shown that $\nu \in \mathcal{M}$ is $\boxplus$-infinitely divisible if and only if $\phi_\nu$ has an analytic continuation to $\mathbf{C}^+$ with values in $\mathbf{C}^- \cup \mathbf{R}$. The Nevanlinna representation of functions with negative imaginary part yields then the following free analogue of the Lévy-Hinčin formula.

3.2. THEOREM. *A measure $\nu \in \mathcal{M}$ is $\boxplus$-infinitely divisible if and only if there exist a finite positive Borel measure $\sigma$ on $\mathbf{R}$, and a real number $\gamma$ such that*
$$\phi_\nu(z) = \gamma + \int_{-\infty}^\infty \frac{1+tz}{z-t}\, d\sigma(t), \quad z \in \mathbf{C}^+.$$

The analogy between Theorems 3.1 and 3.2 is rather formal, though in special cases it was known to go deeper. For instance, if $\sigma$ is a Dirac mass at zero, the formula in Theorem 3.1 yields a Gaussian measure, while the formula in Theorem 3.2 gives a semicircle law. The Gaussian is related with the classical central limit theorem, while the semicircle is related with the free central limit theorem as first shown by Voiculescu [17]. The purpose of this section is to show that this analogy extends to all infinitely divisible measures. In order to facilitate the statements of our results we introduce some additional notation. Fix a finite positive Borel measure $\sigma$ on $\mathbf{R}$ and a real number $\gamma$. We denote by $\nu_*^{\gamma,\sigma}$ the $*$-infinitely divisible measure determined by the formula
$$(\mathcal{F}\nu_*^{\gamma,\sigma})(t) = \exp\left[i\gamma t + \int_{-\infty}^\infty (e^{itx} - 1 - itx)\frac{x^2+1}{x^2}\, d\sigma(x)\right], \quad t \in \mathbf{R},$$
and we denote by $\nu_\boxplus^{\gamma,\sigma}$ the $\boxplus$-infinitely divisible measure such that
$$\phi_{\nu_\boxplus^{\gamma,\sigma}}(z) = \gamma + \int_{-\infty}^\infty \frac{1+tz}{z-t}\, d\sigma(t), \quad z \in \mathbf{C}^+.$$

The following result is essentially contained in [8, §§18 and 24].

3.3. THEOREM. *Let $\nu_*^{\gamma,\sigma}$ be a $*$-infinitely divisible measure, let $\mu_n \in \mathcal{M}$ be a sequence, and let $k_1 < k_2 < \cdots$ be natural numbers. The following assertions are equivalent*:

(1) *the sequence $\underbrace{\mu_n * \mu_n * \cdots * \mu_n}_{k_n \text{ times}}$ converges weakly to $\nu_*^{\gamma,\sigma}$;*

(2) *the measures*
$$d\sigma_n(x) = k_n \frac{x^2}{x^2+1}\, d\mu_n(x)$$
*converge weakly to $\sigma$ and*
$$\lim_{n\to\infty} k_n \int_{-\infty}^\infty \frac{x}{x^2+1}\, d\mu_n(x) = \gamma.$$

The preceding theorem is not stated as such in [8] because the results there are concerned mostly with convolutions of distinct measures. It is interesting to note that the numbers $k_n \int_{-\tau}^{\tau} x \, d\mu_n(x)$ are preferred in [8] to $k_n \int_{-\infty}^{\infty} x/(x^2+1) \, d\mu_n(x)$, most likely on account of their probabilistic interpretation as truncated expected values. It is easy to see that under the equivalent conditions of Theorem 3.3 we have

$$k_n \int_{-\tau}^{\tau} x \, d\mu_n(x) - k_n \int_{-\infty}^{\infty} \frac{x}{x^2+1} \, d\mu_n(x) \underset{n\to\infty}{\longrightarrow} \int_{-\tau}^{\tau} x \, d\sigma(x) - \int_{|x|>\tau} \frac{1}{x} \, d\sigma(x),$$

provided that $\tau$ and $-\tau$ are not atoms of $\sigma$.

We are now ready to state the main result of this section.

3.4. THEOREM. *Fix a finite positive Borel measure $\sigma$ on $\mathbf{R}$, a real number $\gamma$, a sequence $\mu_n \in \mathcal{M}$, and a sequence of positive integers $k_1 < k_2 < \cdots$. The following assertions are equivalent*:

(1) *the sequence* $\underbrace{\mu_n * \mu_n * \cdots * \mu_n}_{k_n \text{ times}}$ *converges weakly to $\nu_*^{\gamma,\sigma}$*;

(2) *the sequence* $\underbrace{\mu_n \boxplus \mu_n \boxplus \cdots \boxplus \mu_n}_{k_n \text{ times}}$ *converges weakly to $\nu_\boxplus^{\gamma,\sigma}$*;

(3) *the measures*

$$d\sigma_n(x) = k_n \frac{x^2}{x^2+1} \, d\mu_n(x)$$

*converge weakly to $\sigma$ and*

$$\lim_{n\to\infty} k_n \int_{-\infty}^{\infty} \frac{x}{x^2+1} \, d\mu_n(x) = \gamma.$$

*Proof.* The equivalence of (1) and (3) is just a restatement of Theorem 3.3. We will prove the equivalence of (2) and (3). Assume first that (2) holds. By Lemma 2.2 and Proposition 2.3 there exists a truncated cone $\Gamma_{\eta,M}$ such that

$$\lim_{n\to\infty} k_n \phi_{\mu_n}(z) = \phi_{\nu_\boxplus^{\gamma,\sigma}}(z), \quad z \in \Gamma_{\eta,M},$$

and $|k_n \phi_{\mu_n}(z)| \leq u(z)$ with $\lim_{z\to\infty, z\in\Gamma_{\eta,M}} u(z)/z = 0$. Then we have

$$|\phi_{\mu_n}(z)| \leq \frac{u(z)}{k_n} \underset{n\to\infty}{\longrightarrow} 0,$$

and therefore $\mu_n$ converges weakly to $\delta_0$; indeed, $\phi_{\delta_0} = 0$. From Proposition 2.6 we get

$$\phi_{\mu_n}(z) = z^2 \left[ G_{\mu_n}(z) - \frac{1}{z} \right] (1 + v_n(z)),$$

where $|v_n(z)| \leq v(z)$, and $\lim_{z \to \infty, z \in \Gamma_{\eta,M}} v(z) = 0$. By decreasing the cone we may assume that $v(z) < 1/2$, and using Proposition 2.7 we deduce that

$$(3.5) \quad k_n z^2 \left[ G_{\mu_n}(z) - \frac{1}{z} \right] = \frac{k_n \phi_{\mu_n}(z)}{1 + v_n(z)} \xrightarrow[n \to \infty]{} \phi_{\nu_{\boxplus}^{\gamma,\sigma}}(z), \quad z \in \Gamma_{\eta,M},$$

and

$$(3.6) \quad \left| k_n z^2 \left[ G_{\mu_n}(z) - \frac{1}{z} \right] \right| \leq \frac{u(z)}{1 + v(z)} \leq 2u(z), \quad z \in \Gamma_{\eta,M}.$$

Observe now that

$$k_n z^2 \left[ G_{\mu_n}(z) - \frac{1}{z} \right] = \int_{-\infty}^{\infty} \frac{k_n t z}{z - t} \, d\mu_n(t), \quad z \in \mathbf{C}^+.$$

It is easily seen that this function has negative imaginary part, though it is not written in standard Nevanlinna form. We can write it in that form as follows:

$$k_n z^2 \left[ G_{\mu_n}(z) - \frac{1}{z} \right] = \gamma_n + \int_{-\infty}^{\infty} \frac{1 + tz}{z - t} \, d\sigma_n(t), \quad z \in \mathbf{C}^+,$$

where

$$\gamma_n = k_n \int_{-\infty}^{\infty} \frac{t}{1 + t^2} \, d\mu_n(t),$$

and

$$d\sigma_n(t) = \frac{k_n t^2}{1 + t^2} \, d\mu_n(t).$$

Inequality (3.6) for $z = iy$, $y \in \mathbf{R}$, implies now

$$\frac{1}{2} \sigma_n(\{t : |t| \geq y\}) \leq \int_{-\infty}^{\infty} \frac{1 + t^2}{y^2 + t^2} \, d\sigma_n(t)$$
$$= -\frac{k_n}{y} \Im \left( z^2 \left[ G_{\mu_n}(z) - \frac{1}{z} \right] \right)$$
$$\leq \frac{2u(iy)}{y}.$$

Because $\lim_{y \to \infty} u(iy) = 0$ we conclude that the sequence $\sigma_n$ is tight, and hence it has weakly convergent subsequences. If a subsequence $\sigma_{n_j}$ converges weakly to a measure $\sigma'$ and $z = x + iy \in \Gamma_{\eta,M}$, then from (3.5) we must have

$$\int_{-\infty}^{\infty} \frac{y}{y^2 + (x - t)^2} (1 + t^2) \, d\sigma'(t) = \lim_{j \to \infty} \int_{-\infty}^{\infty} \frac{y}{y^2 + (x - t)^2} (1 + t^2) \, d\sigma_{n_j}(t)$$
$$= -\lim_{j \to \infty} k_{n_j} \Im \left( z^2 \left[ G_{\mu_n}(z) - \frac{1}{z} \right] \right)$$
$$= -\Im \phi_{\nu_{\boxplus}^{\gamma,\sigma}}(iy)$$
$$= \int_{-\infty}^{\infty} \frac{y}{y^2 + (x - t)^2} (1 + t^2) \, d\sigma(t).$$

Thus the measures $(1 + t^2) \, d\sigma'(t)$ and $(1 + t^2) \, d\sigma(t)$ have the same Poisson integral (in an open set and hence everywhere in $\mathbf{C}^+$), and this can only happen if $\sigma' = \sigma$. We deduce that the sequence $\sigma_n$ converges weakly to $\sigma$. Finally,

$$\gamma = \phi_{\nu_\boxplus^{\gamma,\sigma}}(z) - \int_{-\infty}^{\infty} \frac{1+tz}{z-t} \, d\sigma(t)$$

$$= \lim_{n \to \infty} \left\{ \phi_{\nu_\boxplus^{\gamma,\sigma}}(z) - \int_{-\infty}^{\infty} \frac{1+tz}{z-t} \, d\sigma_n(t) \right\}$$

$$= \lim_{n \to \infty} \left\{ \gamma_n + \phi_{\nu_\boxplus^{\gamma,\sigma}}(z) - k_n z^2 \left[ G_\mu(z) - \frac{1}{z} \right] \right\}.$$

Since

$$\phi_{\nu_\boxplus^{\gamma,\sigma}}(z) - k_n z^2 \left[ G_{\mu_n}(z) - \frac{1}{z} \right] \xrightarrow[n \to \infty]{} 0,$$

we conclude that $\lim_{n \to \infty} \gamma_n$ exists and it equals $\gamma$.

Conversely, assume that (3) holds. Since

$$\mu_n(\{t : |t| > \varepsilon\}) \leq \frac{1+\varepsilon^2}{\varepsilon^2} \int_{-\infty}^{\infty} \frac{t^2}{1+t^2} \, d\mu_n(t) \leq \frac{1+\varepsilon^2}{\varepsilon^2} \frac{1}{k_n} \sigma_n(\mathbf{R})$$

for every $\varepsilon > 0$, we deduce that $\mu_n$ converges weakly to $\delta_0$. Therefore there exists a truncated cone $\Gamma_{\eta,M}$ such that

$$\phi_{\mu_n}(z) = z^2 \left[ G_{\mu_n}(z) - \frac{1}{z} \right] (1 + v_n(z)), \quad z \in \Gamma_{\eta,M},$$

where the $v_n$ are dominated by a positive function $v$ with $\lim_{z \to \infty, z \in \Gamma_{\eta,M}} v(z) = 0$. Using again the Nevanlinna representation

$$k_n z^2 \left[ G_{\mu_n}(z) - \frac{1}{z} \right] = \gamma_n + \int_{-\infty}^{\infty} \frac{1+tz}{z-t} \, d\sigma_n(t), \quad z \in \mathbf{C}^+,$$

we see that condition (3) implies immediately that

$$\lim_{n \to \infty} k_n z^2 \left[ G_{\mu_n}(z) - \frac{1}{z} \right] = \phi_{\nu_\boxplus^{\gamma,\sigma}}(z), \quad z \in \mathbf{C}^+.$$

Hence

$$k_n \phi_{\mu_n}(z) = k_n z^2 \left[ G_{\mu_n}(z) - \frac{1}{z} \right] (1 + v_n(z))$$

also converges to $\phi_{\nu_\boxplus^{\gamma,\sigma}}(z)$ as $n \to \infty$ for $z \in \Gamma_{\eta,M}$. To conclude the proof it suffices (by Prop. 2.3) to show that

$$k_n(iy)^2 \left[ G_{\mu_n}(iy) - \frac{1}{iy} \right] = o(y) \quad \text{as } y \to \infty$$

uniformly in $n$. Since $\gamma_n$ is a bounded sequence, this amounts to showing that

$$\int_{-\infty}^{\infty} \frac{1+ity}{iy-t} \, d\sigma_n(t) = o(y) \quad \text{as } y \to \infty$$

uniformly in $n$. It is easy to verify that
$$\left|\frac{1+ity}{iy-t}\right| \leq y$$
for $y \geq 1$. Hence, if $M > 0$ and $y \geq 1$,
$$\left|\int_{-\infty}^{\infty} \frac{1+ity}{iy-t} d\sigma_n(t)\right| \leq \int_{-M}^{M} 2\frac{1+|t|y}{y+|t|} d\sigma_n(t) + y\sigma_n(\{t : |t| \geq M\})$$
$$\leq 2\frac{1+My}{y+M} + y\sigma_n(\{t : |t| \geq M\}).$$

The desired conclusion follows from the last inequality since $\sigma_n$ was assumed to be a tight sequence. $\square$

## 4. Stable laws and domains of attraction

In this section we will prove the results announced in the introduction on (partial) domains of attraction. We will also discuss briefly the domains of normal attraction of the stable laws. We begin by applying the main result of Section 3 to identically distributed random variables.

4.1. THEOREM. *Suppose we are given a measure $\mu \in \mathcal{M}$, sequences of independent random variables $X_n$ with distribution $\mu$, free random variables $Y_n$ with distribution $\mu$, positive numbers $B_n$, real numbers $A_n$, and natural numbers $k_1 < k_2 < \cdots$. The following two statements are equivalent:*

(1) *the variables*
$$Z_n = \frac{X_1 + X_2 + \cdots + X_{k_n}}{B_n} - A_n$$
*converge in distribution;*

(2) *the variables*
$$W_n = \frac{Y_1 + Y_2 + \cdots + Y_{k_n}}{B_n} - A_n$$
*converge in distribution.*

*If these two equivalent statements are true then there exist a real number $\gamma$, and a positive, finite Borel measure $\sigma$ on $\mathbf{R}$, such that the limit of $Z_n$ is $\nu_*^{\gamma,\sigma}$ and the limit of $W_n$ is $\nu_{\boxplus}^{\gamma,\sigma}$.*

*Proof.* Define measures $\mu_n \in \mathcal{M}$ by
$$\mu_n(S) = \mu\left(B_n S + \frac{A_n}{k_n}\right)$$

for every Borel set $S \subset \mathbf{R}$. Since the distribution of $Z_n$ is

$$\underbrace{\mu_n * \mu_n * \cdots * \mu_n}_{k_n \text{ times}},$$

while the distribution of $W_n$ is

$$\underbrace{\mu_n \boxplus \mu_n \boxplus \cdots \boxplus \mu_n}_{k_n \text{ times}},$$

the result easily follows from Theorem 3.4. □

The preceding result does imply the classical result of Hinčin that a measure $\nu$ is $*$-infinitely divisible if and only if it has a (nonempty) $*$-domain of partial attraction. It also implies the free version of this result, first proved in [13]. Theorem 4.1 gives an identification of the domains of partial attraction as follows.

4.2. COROLLARY. *For every real number $\gamma$ and every finite, positive Borel measure $\sigma$ on $\mathbf{R}$ we have*

$$\mathcal{P}_*(\nu_*^{\gamma,\sigma}) = \mathcal{P}_\boxplus(\nu_\boxplus^{\gamma,\sigma}).$$

Taking $k_n = n$ in Theorem 4.1 we obtain the following result.

4.3. COROLLARY. *Let $\gamma$ be a real number, and $\sigma$ a finite, positive Borel measure on $\mathbf{R}$. The measure $\nu_*^{\gamma,\sigma}$ is $*$-stable if and only if the measure $\nu_\boxplus^{\gamma,\sigma}$ is $\boxplus$-stable. We always have*

$$\mathcal{D}_*(\nu_*^{\gamma,\sigma}) = \mathcal{D}_\boxplus(\nu_\boxplus^{\gamma,\sigma}).$$

Since a stable law always belongs to its own domain of attraction, we have the following "practical" way of determining the correspondence between $*$-stable and $\boxplus$-stable laws.

4.4. COROLLARY. *Assume that $\nu_*^{\gamma,\sigma}$ is $*$-stable (or $\nu_\boxplus^{\gamma,\sigma}$ is $\boxplus$-stable). Then we have that $\nu_*^{\gamma,\sigma} \in \mathcal{D}_\boxplus(\nu_\boxplus^{\gamma,\sigma})$ and $\nu_\boxplus^{\gamma,\sigma} \in \mathcal{D}_*(\nu_*^{\gamma,\sigma})$.*

Given $\nu \in \mathcal{M}$ we denote its dilation by a factor $a > 0$ by $D_a\nu(S) = \nu(aS)$. Let us fix now a $*$-stable measure $\nu$. Associated with $\nu$ there is a number $\alpha \in (0, 2]$ such that the measure $\nu * \nu$ is a translate of the measure $D_{1/2^\alpha}\nu$. The number $\alpha$ is called the *stability index* of $\nu$. Let now $\mu$ belong to the $*$-domain of attraction of $\nu$. Thus, if $X_n$ is a sequence of independent identically distributed variables with common distribution $\mu$, there exist constants $A_n, B_n$, with $B_n > 0$, such that the variables $B_n^{-1}(X_1 + X_2 + \cdots + X_n) - A_n$ converge to $\nu$ in distribution. If the constants $B_n$ can be chosen of the form $cn^{1/\alpha}$, with $c$ independent of $n$, then $\mu$ is in the *normal $*$-domain of attraction* of $\nu$. The normal $\boxplus$-domain of attraction of a $\boxplus$-stable law is defined in an analogous manner, and the following result is then an immediate consequence of Theorem 4.1.

4.5. COROLLARY. *Assume that $\nu_*^{\gamma,\sigma}$ is $*$-stable. Then the normal $*$-domain of attraction of $\nu_*^{\gamma,\sigma}$ coincides with the normal $\boxplus$-domain of attraction of $\nu_\boxplus^{\gamma,\sigma}$.*

The case $\alpha = 2$ of this result was proved in [14], where it was shown that the normal $\boxplus$-domain of attraction of the semicircle laws coincides with the class of measures with finite variance.

## 5. Cauchy transforms and domains of attraction

In order to discuss a little further the ($*$ or $\boxplus$)-domains of attraction we need a few facts from the theory of regularly varying functions. Since these facts, though classical, may not be so well known, we list them here.

A positive Borel function $\ell$ defined in a neighborhood of $+\infty$ is said to be *slowly varying* if
$$\lim_{x \to \infty} \frac{\ell(tx)}{\ell(x)} = 1, \quad \text{for all } t > 0.$$
A positive Borel function $r$ defined in a neighborhood of $+\infty$ is *regularly varying* with *index* $\rho \in \mathbf{R}$ if $r(y) = y^\rho \ell(y)$, with $\ell$ slowly varying. We denote by $\mathcal{R}_\rho$ the class of (germs at $+\infty$ of) regularly varying functions with index $\rho$, so that $\mathcal{R}_0$ is the class of slowly varying functions. As usual, we write $f \sim g$ as $x \to \infty$ if $\lim_{x \to \infty} f(x)/g(x) = 1$. It will be convenient to use the following convention: the asymptotic relation $f \sim kg$ means $\lim_{x \to \infty} f(x)/g(x) = 0$ if $k = 0$. This is useful in writing certain statements in a uniform way.

The following result is due to Karamata. It can be found in [6] (cf. Theorems 1.6.4 and 1.6.5 and the comments following 1.6.5).

5.1. THEOREM. *Let $\alpha$ and $n$ be real numbers such that $0 < \alpha < n$, and let $\rho$ be a finite positive Borel measure on $[0, +\infty)$. The following assertions are equivalent*:
(1) *the function $y \mapsto \rho((y, \infty))$ varies regularly with index $-\alpha$;*
(2) *the function $y \mapsto \int_0^y t^n \, d\rho(t)$ varies regularly with index $n - \alpha$.*
*If these equivalent conditions are verified, then*
$$\int_0^y t^n \, d\rho(t) \sim \frac{\alpha}{n - \alpha} y^n \rho((y, \infty)) \quad as \ y \to \infty.$$

The situation is more complicated when $n = \alpha$, and it was studied by de Haan (see [6, Chap. 3]). We only state a simple result for $\alpha = 1$ because it is the only case used here.

5.2. PROPOSITION. *Let $\rho$ be a finite Borel measure on $[0,+\infty)$ such that the function $y \mapsto \rho((y,\infty))$ varies regularly with index $-1$. Then the function $y \mapsto \int_0^y t\,d\rho(t)$ is slowly varying and*

$$y\rho((y,\infty)) = o\left(\int_0^y t\,d\rho(t)\right) \quad \text{as } y \to \infty.$$

Another useful result of Karamata refers to Stieltjes transforms (cf. [6, Th. 1.7.4]). We recall that a Baire measure is a Borel measure which takes finite values on compact sets.

5.3. THEOREM. *Let $\rho$ be a positive Baire measure on $[0,+\infty)$, and fix $\alpha \in (0,1]$. The following conditions are equivalent:*

(1) *the function $y \mapsto \rho([0,y))$ varies regularly with index $1-\alpha$;*
(2) *the function $y \mapsto \int_0^\infty 1/(t+y)\,d\rho(t)$ varies regularly with index $-\alpha$.*

*If these equivalent conditions are satisfied then*

$$\int_0^\infty \frac{1}{t+y}\,d\rho(t) \sim \frac{(1-\alpha)\pi}{\sin \pi\alpha} \frac{\rho([0,y))}{y} \quad \text{as } y \to \infty,$$

*where the constant $(1-\alpha)\pi/\sin \pi\alpha$ must be replaced by 1 if $\alpha = 1$.*

We actually need a variation of this result.

5.4. COROLLARY. *Let $\rho$ be a positive Baire measure on $[0,+\infty)$, and fix $\alpha \in (0,2]$. The following conditions are equivalent:*

(1) *the function $y \mapsto \rho([0,y))$ varies regularly with index $2-\alpha$;*
(2) *the function $y \mapsto \int_0^\infty 1/(t^2+y^2)\,d\rho(t)$ varies regularly with index $-\alpha$.*

*If these equivalent conditions are satisfied then*

$$\int_0^\infty \frac{1}{t^2+y^2}\,d\rho(t) \sim \frac{(1-\frac{\alpha}{2})\pi}{\sin \frac{\pi\alpha}{2}} \frac{\rho([0,y))}{y^2} \quad \text{as } y \to \infty,$$

*where the constant $(1-\alpha/2)\pi/\sin(\pi\alpha/2)$ must be replaced by 1 if $\alpha = 2$.*

The corollary is readily obtained by observing that

$$\int_0^\infty \frac{1}{t^2+y^2}\,d\rho(t) = \int_0^\infty \frac{1}{t+y_1}\,d\rho_1(t),$$

with $y_1 = y^2$ and $d\rho_1(t) = d\rho(\sqrt{t})$. The proof of the following result is equally easy.

5.5. COROLLARY. *Let $\rho$ and $\alpha$ be as in Corollary 5.4. We have*

$$y^{\alpha-2}\rho([0,y)) = o(1) \quad \text{as } y \to \infty$$

*if and only if*

$$y^\alpha \int_0^\infty \frac{1}{t^2+y^2} \, d\rho(t) = o(1) \quad \text{as } y \to \infty.$$

Aside from the $*$-domain of attraction of the normal law, $*$-domains of attraction (and hence $\boxplus$-domains of attraction) are classified by two parameters $\alpha \in (0,2)$, i.e., the stability index, and $\theta \in [-1,1]$. For further reference we record the following.

5.6. *Definition.* Fix $\alpha \in (0,2)$ and $\theta \in [-1,1]$. We say that a measure $\mu \in \mathcal{M}$ belongs to the domain $\mathcal{C}_{\alpha,\theta}$ if the following conditions are satisfied:

(1) the function $y \mapsto \int_{-y}^y t^2 \, d\mu(t)$ varies regularly with index $2-\alpha$;
(2)
$$\lim_{t\to\infty} \frac{\mu((t,\infty)) - \mu((-\infty,-t))}{\mu((t,\infty)) + \mu((-\infty,-t))} = \theta.$$

Observe that for $\mu \in \mathcal{C}_{\alpha,\theta}$,

$$\lim_{t\to\infty} \frac{\mu((t,\infty))}{\mu((t,\infty)) + \mu((-\infty,-t))} = p,$$

$$\lim_{t\to\infty} \frac{\mu((-\infty,-t))}{\mu((t,\infty)) + \mu((-\infty,-t))} = q,$$

where $p = (1+\theta)/2$ and $q = (1-\theta)/2$ so that $p \geq 0$, $q \geq 0$, $p+q = 1$, and $p - q = \theta$.

A classical result on domains of attraction (see [8]) is that each $\mathcal{C}_{\alpha,\theta}$ contains exactly one equivalence class of $*$-stable measures, and $\mathcal{C}_{\alpha,\theta}$ coincides with the $*$-domain of attraction of each of those $*$-stable measures. By virtue of the results of Section 4 we have then proved the following theorem.

5.7. THEOREM. *Let $\alpha \in (0,2)$ and $\theta \in [-1,1]$. Then there exist a unique class of equivalence of $*$-stable laws $[\nu]$ and a unique class of equivalence of $\boxplus$-stable laws $[\nu']$ such that $\nu \in \mathcal{C}_{\alpha,\theta}$ and $\nu' \in \mathcal{C}_{\alpha,\theta}$. Moreover $\mathcal{D}_*(\nu) = \mathcal{D}_\boxplus(\nu') = \mathcal{C}_{\alpha,\theta}$.*

The above theorem does not take into account the case corresponding to $\alpha = 2$, i.e., the central limit theorem (see [14]). Indeed a measure $\mu \in \mathcal{M}$ belongs to the $*$-domain of attraction of the Gaussian if and only if it belongs to the $\boxplus$-domain of attraction of the semicircle law if and only if the function $y \mapsto \int_{-y}^y t^2 \, d\mu(t)$ is slowly varying.

We will describe measures in $\mathcal{C}_{\alpha,\theta}$ via the values of their Cauchy transforms on the imaginary axis. We begin by observing that

$$G_\mu(iy) - \frac{1}{iy} = \frac{i}{y}\int_{-\infty}^\infty \frac{t^2}{y^2+t^2}\,d\mu(t) - \int_{-\infty}^\infty \frac{t}{y^2+t^2}\,d\mu(t)$$

for $\mu \in \mathcal{M}$ and $y > 0$.

5.8. PROPOSITION. *For a measure $\mu \in \mathcal{M}$ and $\alpha \in (0,2)$ the following assertions are equivalent*:

(1) *the function $y \mapsto \mu(\{t : |t| > y\})$ varies regularly with index $-\alpha$*;
(2) *the function $y \mapsto \int_{-y}^y t^2\,d\mu(t)$ varies regularly with index $2-\alpha$*;
(3) *the function $y \mapsto \int_{-y}^y |t|^3\,d\mu(t)$ varies regularly with index $3-\alpha$*;
(4) *the function $y \mapsto \int_{-\infty}^\infty t^2/(y^2+t^2)\,d\mu(t)$ varies regularly with index $-\alpha$*.

*If $\alpha \in (1,2)$ then (1)–(4) are also equivalent to*:

(5) *the function $y \mapsto \int_{-\infty}^\infty |t|^3/(y^2+t^2)\,d\mu(t)$ varies regularly with index $1-\alpha$*.

*If $\alpha \in (0,1)$ then conditions (1)–(4) are also equivalent to*:

(6) *the function $y \mapsto \int_{-y}^y |t|\,d\mu(t)$ varies regularly with index $1-\alpha$*;
(7) *the function $y \mapsto \int_{-\infty}^\infty |t|^3/(y^2+t^2)\,d\mu(t)$ varies regularly with index $-1-\alpha$*.

*If $\alpha = 1$ then condition (1) implies conditions (6) and (7), which are equivalent to each other. If (1) is satisfied, we have*

$$\int_{-\infty}^\infty \frac{t^2}{y^2+t^2}\,d\mu(t) \sim \frac{(1-\frac{\alpha}{2})\pi}{\sin\frac{\pi\alpha}{2}}\frac{1}{y^2}\int_{-y}^y t^2\,d\mu(t)$$
$$\sim \frac{\frac{\pi\alpha}{2}}{\sin\frac{\pi\alpha}{2}}\mu(\{t:|t|>y\}) \quad \text{as } y\to\infty,$$

*and*

$$\frac{1}{y^3}\int_{-y}^y |t|^3\,d\mu(t) \sim \frac{\alpha}{3-\alpha}\mu(\{t:|t|>y\}) \quad \text{as } y\to\infty.$$

*If $\alpha \in (1,2)$ we also have*

$$\int_{-\infty}^\infty \frac{|t|^3}{y^2+t^2}\,d\mu(t) \sim \frac{\frac{3-\alpha}{2}\pi}{-\cos\frac{\pi\alpha}{2}}\frac{1}{y^2}\int_{-y}^y |t|^3\,d\mu(t)$$
$$\sim \frac{\frac{\pi\alpha}{2}}{-\cos\frac{\pi\alpha}{2}}y\mu(\{t:|t|>y\}) \quad \text{as } y\to\infty.$$

*If $\alpha \in (0,1)$ then*

$$\int_{-\infty}^\infty \frac{|t|}{y^2+t^2}\,d\mu(t) \sim \frac{\frac{1-\alpha}{2}\pi}{\cos\frac{\pi\alpha}{2}}\frac{1}{y^2}\int_{-y}^y |t|\,d\mu(t)$$
$$\sim \frac{\frac{\pi\alpha}{2}}{\cos\frac{\pi\alpha}{2}}\frac{1}{y}\mu(\{t:|t|>y\}) \quad \text{as } y\to\infty.$$

*Finally, if $\alpha = 1$,*

$$\int_{-\infty}^{\infty} \frac{|t|}{y^2 + t^2} \, d\mu(t) \sim \frac{1}{y^2} \int_{-y}^{y} |t| \, d\mu(t) \quad \text{as } y \to \infty,$$

*and*

$$\mu(\{t : |t| > y\}) = o\left(\frac{1}{y} \int_{-y}^{y} |t| \, d\mu(t)\right) \quad \text{as } y \to \infty.$$

*Proof.* The equivalence of (1), (2), (3), and (6) (for $\alpha \in (0,1)$) follows immediately from Theorem 5.1 applied to the measure $\rho$ defined by

$$\rho(S) = \mu(\{t : |t| \in S\}), \quad S \subset [0, +\infty).$$

Proposition 5.2 applied to the same measure shows that (1) implies (6) if $\alpha = 1$. In order to prove the equivalence of (3) and (4) consider the Baire measure $\rho$ on $[0, \infty)$ determined by

$$\rho([0, y)) = \int_{(-y, y)} t^2 \, d\mu(t), \quad y > 0.$$

Then we can write

$$\int_{-\infty}^{\infty} \frac{t^2}{y^2 + t^2} \, d\mu(t) = \int_{0}^{\infty} \frac{d\rho(t)}{y^2 + t^2},$$

and the desired equivalence follows immediately from Corollary 5.4. The equivalence between (3) and (5) (for $\alpha \in (1,2)$) and between (6) and (7) (for $\alpha \in (0,1]$) is proved analogously. The asymptotic relations follow from the corresponding results about slowly varying functions. □

We can now characterize measures in $\mathcal{C}_{\alpha,\theta}$ in terms of the behavior of their Cauchy transforms on the imaginary axis. We prove three statements corresponding to the cases $\alpha < 1$, $\alpha > 1$, and $\alpha = 1$.

5.9. PROPOSITION. *Fix $\alpha \in (0,1)$, $\theta \in [-1,1]$, and $\mu \in \mathcal{M}$. The following assertions are equivalent:*

(1) $\mu \in \mathcal{C}_{\alpha,\theta}$;
(2) *there exists $f \in \mathcal{R}_{-\alpha-1}$ such that*

$$G_\mu(iy) - \frac{1}{iy} = \left(i + \theta \tan \frac{\pi\alpha}{2}\right) f(y)(1 + o(1)) \quad \text{as } y \to \infty;$$

(3) *there exists $g \in \mathcal{R}_{1-\alpha}$ such that*

$$\phi_\mu(iy) = -\left(i + \theta \tan \frac{\pi\alpha}{2}\right) g(y)(1 + o(1)) \quad \text{as } y \to \infty.$$

*Proof.* We observe first that the equivalence between (2) and (3) is an immediate consequence of Proposition 2.5; in fact, we can take $g(y) = y^2 f(y)$. We must therefore prove the equivalence of (1) and (2). Assume first that (1) holds, and set
$$f(y) = \frac{1}{y}\int_{-\infty}^{\infty} \frac{t^2}{y^2+t^2}\,d\mu(t), \quad y > 0.$$
The function $f$ belongs to $\mathcal{R}_{-1-\alpha}$ by Proposition 5.8, and
$$\mu(\{t : |t| > y\}) \sim \frac{\sin\frac{\pi\alpha}{2}}{\frac{\pi\alpha}{2}} y f(y) \quad \text{as } y \to \infty.$$
Furthermore, since $\alpha \in (0,1)$, we also have
$$\int_{-\infty}^{\infty} \frac{|t|}{y^2+t^2}\,d\mu(t) \sim \frac{\frac{\pi\alpha}{2}}{\cos\frac{\pi\alpha}{2}}\frac{1}{y}\mu(\{t : |t| > y\})$$
$$\sim \tan\frac{\pi\alpha}{2} f(y) \quad \text{as } y \to \infty.$$
Consider now the numbers $p = (1+\theta)/2$ and $q = (1-\theta)/2$. Since $\mu \in \mathcal{C}_{\alpha,\theta}$,
$$\mu((y,\infty)) \sim p\mu(\{t : |t| > y\}) \quad \text{as } y \to \infty.$$
Thus $\mu((y,\infty))$ varies regularly with index $-\alpha$ (if $p \neq 0$), and
$$\int_0^{\infty} \frac{t}{y^2+t^2}\,d\mu(t) \sim \frac{\frac{\pi\alpha}{2}}{\cos\frac{\pi\alpha}{2}}\frac{1}{y}\mu((y,\infty))$$
$$\sim p\frac{\frac{\pi\alpha}{2}}{\cos\frac{\pi\alpha}{2}}\frac{1}{y}\mu(\{t : |t| > y\})$$
$$\sim p\int_{-\infty}^{\infty} \frac{|t|}{y^2+t^2}\,d\mu(t) \quad \text{as } y \to \infty.$$
If $p = 0$ this relation simply means
$$\int_0^{\infty} \frac{t}{y^2+t^2}\,d\mu(t) = o\left(\int_{-\infty}^{\infty} \frac{|t|}{y^2+t^2}\,d\mu(t)\right) \quad \text{as } y \to \infty.$$
In an analogous manner,
$$\int_{-\infty}^{0} \frac{t}{y^2+t^2}\,d\mu(t) \sim -q\int_{-\infty}^{\infty} \frac{|t|}{y^2+t^2}\,d\mu(t) \quad \text{as } y \to \infty,$$
so that
$$\int_{-\infty}^{\infty} \frac{t}{y^2+t^2}\,d\mu(t) \sim (p-q)\int_{-\infty}^{\infty} \frac{|t|}{y^2+t^2}\,d\mu(t)$$
$$\sim \theta\tan\frac{\pi\alpha}{2} f(y) \quad \text{as } y \to \infty.$$
The representation in (2) follows immediately from the formula
$$G_\mu(iy) - \frac{1}{iy} = \frac{i}{y}\int_{-\infty}^{\infty} \frac{t^2}{y^2+t^2}\,d\mu(t) - \int_{-\infty}^{\infty} \frac{t}{y^2+t^2}\,d\mu(t).$$

Conversely, assume that (2) holds. Then

$$\Im\left[G_\mu(iy) - \frac{1}{iy}\right] = \frac{i}{y}\int_{-\infty}^{\infty}\frac{t^2}{y^2+t^2}\,d\mu(t)$$

is a function of regular variation with index $-\alpha - 1$, asymptotic to $f(y)$ as $y \to \infty$. Considering real parts we also obtain

$$\int_{-\infty}^{\infty}\frac{t}{y^2+t^2}\,d\mu(t) \sim \theta\tan\frac{\pi\alpha}{2}f(y) \quad \text{as } y \to \infty.$$

Since we know from Proposition 5.8 that

$$\int_{-\infty}^{\infty}\frac{|t|}{y^2+t^2}\,d\mu(t) \sim \tan\frac{\pi\alpha}{2}f(y) \quad \text{as } y \to \infty,$$

we conclude that

$$\int_0^\infty \frac{t}{y^2+t^2}\,d\mu(t) = \frac{1}{2}\int_{-\infty}^\infty \frac{|t|}{y^2+t^2}\,d\mu(t) + \frac{1}{2}\int_{-\infty}^\infty \frac{t}{y^2+t^2}\,d\mu(t)$$
$$\sim \frac{1+\theta}{2}\tan\frac{\pi\alpha}{2}f(y)$$
$$\sim p\frac{\frac{\pi\alpha}{2}}{\cos\frac{\pi\alpha}{2}}\frac{1}{y}\mu(\{t:|t|>y\}) \quad \text{as } y \to \infty.$$

Proposition 5.8 (applied to the restriction of $\mu$ to the positive axis) implies now that

$$\mu((t,\infty)) \sim p\mu(\{t:|t|>y\}) \quad \text{as } y \to \infty,$$

and therefore we can conclude that $\mu \in \mathcal{C}_{\alpha,\theta}$. The proposition is proved. $\square$

If $\alpha > 1$ then measures in $\mathcal{C}_{\alpha,\theta}$ have finite mean, and therefore each measure in that domain is equivalent to a measure with zero mean. Thus there is little loss of generality in considering only measures with zero mean in this case.

5.10. PROPOSITION. *Fix $\alpha \in (1,2)$, $\theta \in [-1,1]$, and $\mu \in \mathcal{M}$. Assume that $\mu$ has zero mean. The following assertions are equivalent*:

(1) $\mu \in \mathcal{C}_{\alpha,\theta}$;
(2) *there exists $f \in \mathcal{R}_{-\alpha-1}$ such that*

$$G_\mu(iy) - \frac{1}{iy} = \left(i + \theta\tan\frac{\pi\alpha}{2}\right)f(y)(1+o(1)) \quad \text{as } y \to \infty;$$

(3) *there exists $g \in \mathcal{R}_{1-\alpha}$ such that*

$$\phi_\mu(iy) = -\left(i + \theta\tan\frac{\pi\alpha}{2}\right)g(y)(1+o(1)) \quad \text{as } y \to \infty.$$

*Proof.* Since $\mu$ has zero mean,

$$\int_{-\infty}^{\infty} \frac{t}{y^2+t^2}\, d\mu(t) = \int_{-\infty}^{\infty} \left[\frac{t}{y^2+t^2} - \frac{t}{y^2}\right] d\mu(t)$$
$$= -\frac{1}{y^2} \int_{-\infty}^{\infty} \frac{t^3}{y^2+t^2}\, d\mu(t).$$

The proof proceeds along the lines of the preceding one, with

$$\int_{-\infty}^{\infty} t/(y^2+t^2)\, d\mu(t)$$

replaced by

$$\int_{-\infty}^{\infty} t^3/(y^2+t^2)\, d\mu(t).$$

We leave the details to the interested reader. □

The situation is more complicated for measures in $\mathcal{C}_{1,\theta}$, some of which have finite mean. We limit ourselves to a partial result.

5.11. PROPOSITION. *Fix $\theta \in [-1, 1]$, and $\mu \in \mathcal{C}_{1,\theta}$. Then*

$$G_\mu(iy) - \frac{1}{iy} = -\frac{1}{y^2} \int_{-y}^{y} t\, d\mu(t) + if(y) + o(1)g(y) \quad \text{as } y \to \infty,$$

*where $f, g \in \mathcal{R}_{-2}$ and satisfy*

$$g(y) \sim \frac{1}{y^2} \int_{-y}^{y} |t|\, d\mu(t) \quad \text{as } y \to \infty$$

*and*

$$f(y) \sim \frac{\pi}{2}\frac{1}{y}\mu(\{t : |t| > y\}) = o(g(y)) \quad \text{as } y \to \infty.$$

*Proof.* As before, $f(y)$ can be taken to be the imaginary part of $G_\mu(iy) - 1/iy$, so we only need to consider the real part. Since

$$y \mapsto \mu((y, +\infty)) \sim p\mu(\{t : |t| > y\}) \quad \text{as } y \to \infty$$

is a function of regular variation with index $-1$ (if $p \neq 0$), Proposition 5.8 implies that

$$\int_0^\infty \frac{t}{y^2+t^2}\, d\mu(t) \sim \frac{1}{y^2} \int_0^y t\, d\mu(t) \quad \text{as } y \to \infty.$$

Thus setting

$$g(y) = \frac{1}{y^2} \int_{-y}^{y} |t|\, d\mu(t) \quad \text{as } y \to \infty,$$

it is easy to conclude that the desired representation holds. □

It is somewhat unsatisfactory that the constant $\theta$ does not appear explicitly in the formula given in Proposition 5.11. Let us note, however, that for measures $\mu$ with infinite mean,

$$\theta y^2 g(y) = \theta \int_{-y}^{y} |t|\, d\mu(t) \sim \int_{-y}^{y} t\, d\mu(t) \quad \text{as } y \to \infty,$$

and the statement takes a form more similar to Propositions 5.9 and 5.10.

It is now fairly easy to determine the $\boxplus$-stable laws in each of the domains $\mathcal{C}_{\alpha,\theta}$. The complex functions in the following statement are given by their principal value in the upper half-plane.

5.12. PROPOSITION.
(1) If $\alpha \in (0,1) \cup (1,2)$ and $\theta \in [-1,1]$, then the measure $\nu \in \mathcal{M}$ with

$$\phi_\nu(z) = -\left(i + \theta \tan \frac{\pi\alpha}{2}\right) i^{\alpha-1} z^{1-\alpha}$$

belongs to $\mathcal{C}_{\alpha,\theta}$.

(2) If $\theta \in [-1,1]$, then the measure $\nu \in \mathcal{M}$ with

$$\phi_\nu(z) = 2\theta \log z - i\pi(1+\theta)$$

belongs to $\mathcal{C}_{1,\theta}$.

*Proof.* For $\alpha \neq 1$ it is known (see [4]) that the Voiculescu transforms of the stable laws of index $\alpha$ have the form $a + bz^{1-\alpha}$. Among these, the function $-i^{\alpha-1} z^{1-\alpha}$ is equal to $y^{1-\alpha}$ for $z = iy$, and this function is of regular variation with index $1 - \alpha$. The general form in the statement is easily deduced from Propositions 5.9 and 5.10. In order to arrive at the formula (2), we note that a measure $\mu$ such that

$$y\mu(\{t : |t| > y\}) \sim 1 \quad \text{as } y \to \infty$$

must satisfy

$$\int_{-y}^{y} |t|\, d\mu(t) \sim \log y \quad \text{as } y \to \infty.$$

The functions in (2) are exactly those stable laws of index 1 which achieve the proper balance between real and imaginary parts. $\square$

If $\mu \in \mathcal{M}$ a $\boxplus$-stable distribution with stability index $\alpha \in (0,2]$, from the relation $\phi_{D_a\mu}(z) = a\phi_\mu(z/a)$ (where $D_a\mu$ is the dilation of $\mu$ by a factor $a > 0$), it is easy to verify that the measure $D_a\mu \boxplus D_b\mu$ is a translate of the measure $D_{(a^\alpha + b^\alpha)^{1/\alpha}}\mu$ (and is equal to it in the case of a strictly $\boxplus$-stable distribution).

It is interesting to note that every ⊞-stable distribution of stability index $\neq 1$ is equivalent to some strictly ⊞-stable distribution, but this is not the case for $\alpha = 1$, where the strictly ⊞-stable distributions are the ones with constant Voiculescu transform, namely, the Dirac measures and the measures equivalent to the Cauchy distribution. Note that these measures are also exactly the strictly ∗-stable distributions of stability index 1 (see e.g. [8]).

## 6. Limit laws for Boolean convolution

We have seen earlier how one can associate to a given measure $\mu \in \mathcal{M}$ its Cauchy transform $G_\mu$ and its reciprocal $F_\mu = 1/G_\mu : \mathbf{C}^+ \to \mathbf{C}^+$. We have $\Im z \leq \Im F_\mu(z)$ so that the function $E_\mu(z) = z - F_\mu(z)$ maps $\mathbf{C}^+$ to $\mathbf{C}^- \cup \mathbf{R}$, and, in addition, $E_\mu(z)/z \to 0$ as $z \to \infty$ nontangentially. Conversely, if $E : \mathbf{C}^+ \to \mathbf{C}^- \cup \mathbf{R}$ is an analytic function so that $E(z)/z \to 0$ as $z \to \infty$ nontangentially, then there exists $\mu \in \mathcal{M}$ such that $E_\mu = E$. This observation leads to the formal definition of the Boolean convolution introduced in [16]. Given $\mu, \nu \in \mathcal{M}$, there exists $\rho \in \mathcal{M}$ such that
$$E_\rho = E_\mu + E_\nu.$$
The measure $\rho$ is called the *Boolean convolution* of $\mu$ and $\nu$, and it is denoted $\mu \uplus \nu$. Boolean convolution is an associative, commutative law, with $\delta_0$ as the zero element. We have $\delta_s \uplus \delta_t = \delta_{s+t}$, but generally $\delta_t \uplus \mu$ is not a translate of $\mu$. We will show that the limit laws of Boolean convolution are determined by the limit laws of classical probability, just as in the case of free convolution. The basic tool is the following counterpart of Propositions 2.5 and 2.6.

6.1. PROPOSITION. *For every $\mu \in \mathcal{M}$,*
$$E_\mu(z) = z^2 \left[ G_\mu(z) - \frac{1}{z} \right] (1 + o(1))$$
*as $z \to \infty$ nontangentially. Moreover, if $\mathcal{F} \subset \mathcal{M}$ is a tight family and we write*
$$E_\mu(z) = z^2 \left[ G_\mu(z) - \frac{1}{z} \right] (1 + v_\mu(z)), \quad z \in \mathbf{C}^+, \mu \in \mathcal{F},$$
*then $\sup_{\mu \in \mathcal{F}} |v_\mu(z)| = o(1)$ as $z \to \infty$ nontangentially.*

*Proof.* Observe that
$$E_\mu(z) = z^2 \left[ G_\mu(z) - \frac{1}{z} \right] \frac{1}{zG_\mu(z)}, \quad z \in \mathbf{C}^+.$$
The conclusion now can be deduced as in the proof of Propositions 2.5 and 2.6, except that the situation is much simpler here because no inversion is required. □

An immediate consequence of Propositions 6.1 and 2.3 is that the map $\mu \mapsto E_\mu$ behaves well relative to weak convergence. We record the result below.

6.2. PROPOSITION. *Let $\mu_n \in \mathcal{M}$ be a sequence. The following assertions are equivalent*:

(i) *$\mu_n$ converges weakly to a probability measure $\mu$;*
(ii) *there exist $\eta, M > 0$ such that the sequence $E_{\mu_n}$ converges uniformly on $\Gamma_{\eta,M}$ to a function $E$, and $E_{\mu_n}(z) = o(|z|)$ uniformly in $n$ as $z \to \infty$, $z \in \Gamma_{\eta,M}$;*
(iii) *There exists $M' > 0$ such that $\lim_{n\to\infty} E_{\mu_n}(iy)$ exists for every $y > M'$, and $E_{\mu_n}(iy) = o(y)$ uniformly in $n$ as $y \to \infty$.*

*Moreover, if* (i) *and* (ii) *are satisfied, we have $E = E_\mu$.*

As we have seen before, analytic functions $E : \mathbf{C}^+ \to \mathbf{C}^- \cup \mathbf{R}$, with $E(z) = o(z)$ as $z \to \infty$ nontangentially, have a Nevanlinna representation

$$E(z) = \gamma + \int_{-\infty}^{\infty} \frac{1+tz}{z-t}\, d\sigma(t), \quad z \in \mathbf{C}^+,$$

where $\gamma$ is a real number and $\sigma$ is a finite, positive Borel measure on $\mathbf{R}$. We will denote by $\nu_{\boxplus}^{\gamma,\sigma}$ the measure in $\mathcal{M}$ determined by

$$E_{\nu_{\boxplus}^{\gamma,\sigma}}(z) = \gamma + \int_{-\infty}^{\infty} \frac{1+tz}{z-t}\, d\sigma(t), \quad z \in \mathbf{C}^+.$$

Every measure $\mu \in \mathcal{M}$ is of the form $\nu_{\boxplus}^{\gamma,\sigma}$. This reflects the fact that all measures in $\mathcal{M}$ are $\boxplus$-infinitely divisible. We are now ready for the main result of this section Its proof is practically identical with that of Theorem 3.4 (whose statement is also repeated here) and is therefore omitted.

6.3. THEOREM. *Fix a finite positive Borel measure $\sigma$ on $\mathbf{R}$, a real number $\gamma$, a sequence $\mu_n \in \mathcal{M}$, and a sequence of positive integers $k_1 < k_2 < \cdots$. The following assertions are equivalent*:

(1) *The sequence $\underbrace{\mu_n * \mu_n * \cdots * \mu_n}_{k_n \text{ times}}$ converges weakly to $\nu_*^{\gamma,\sigma}$;*

(2) *The sequence $\underbrace{\mu_n \boxplus \mu_n \boxplus \cdots \boxplus \mu_n}_{k_n \text{ times}}$ converges weakly to $\nu_{\boxplus}^{\gamma,\sigma}$;*

(3) *The sequence $\underbrace{\mu_n \uplus \mu_n \uplus \cdots \uplus \mu_n}_{k_n \text{ times}}$ converges weakly to $\nu_{\uplus}^{\gamma,\sigma}$;*

(4) *The measures*

$$d\sigma_n(x) = k_n \frac{x^2}{x^2+1} d\mu_n(x)$$

*converge weakly to* $\sigma$, *and*

$$\lim_{n\to\infty} k_n \int_{-\infty}^{\infty} \frac{x}{x^2+1} d\mu_n(x) = \gamma.$$

We will not write out formally the consequences of this result. Note however that it implies that every $\mu \in \mathcal{M}$ has a nonempty ⊞-domain of partial attraction which is identical with the ∗-domain of partial attraction of some ∗-infinitely divisible measure. The laws which have a nonempty ⊞-domain of attraction can also be easily determined. Not all of these however are ⊞-stable. For instance, if $E_\mu(z)$ has the form $a + b \log z$ with $b \neq 0$, then $\mu$ is not ⊞-stable even though it has a ⊞-domain of attraction. The results of [16] concerning the central limit theorem and limit theorems for the Boolean analogues of the Poisson laws also follow from Theorem 6.3.

## Appendix

### The density of free stable distributions

It is easy to see from Proposition 5.12 (see also [4, Thm. 7.5]), that every ⊞-stable distribution is equivalent to a unique distribution whose Voiculescu transform belongs to the following list:

(1) $\phi(z) = z^{-1}$;
(2) $\phi(z) = e^{i(\alpha-2)\rho\pi} z^{-\alpha+1}$ with $1 < \alpha < 2$, $0 \leq \rho \leq 1$;
(3) we consider two subcases
　　(i) $\phi(z) = 0$,
　　(ii) $\phi(z) = -2\rho i + 2(2\rho-1)/\pi \log z$ with $0 \leq \rho \leq 1$;
(4) $\phi(z) = -e^{i\alpha\rho\pi} z^{-\alpha+1}$ with $0 < \alpha < 1$, $0 \leq \rho \leq 1$.

The stability index of a ⊞-stable distribution is equal to 2 in case (1), to $\alpha$ in cases (2) and (4), and to 1 in case (3). The parameter $\rho$ which appears in cases (2), (3) and (4) will be called the *asymmetry coefficient*, and one can see that the measure corresponding to the parameters $(\alpha, \rho)$ is the image of the measure with parameters $(\alpha, 1-\rho)$ by the map $t \mapsto -t$ on $\mathbf{R}$.

The parameter $\rho$ is related to the parameter $\theta$ considered earlier in this paper by the formula $2\rho = \theta + 1$. In the sequel we will consider only the ⊞-stable distributions whose Voiculescu transforms belong to the list above, and we will denote the ⊞-stable distribution with stability index $\alpha$ and asymmetry coefficient $\rho$ by $\nu_{\alpha,\rho}$.

In some special cases, it is possible to give a simple formula for the ⊞-stable distributions, namely, the case (1) corresponds to the semicircle distribution, the case (3)(i) to the Dirac measure at zero, the case (3)(ii), with $\rho = 1/2$, to the Cauchy distribution, and the density in the cases $(1/2, 0)$ and $(1/2, 1)$ can also be computed in closed form (see the remark after Prop. A1.4).

In Section A1 below we will give a formula for the density of the ⊞-stable distributions from which one can answer almost any question on these distributions. In particular, we will compute explicitly the support of these distributions, and show that they have analytic densities. We will compute the asymptotic behavior of these densities in Section A2, and we will also show that the ⊞-stable distributions are unimodal. In Section A3, we will see that they satisfy a duality relation analogous to that of the classical stable distributions on the real line (we refer to the monograph [23] of Zolotarev for information on the stable distributions on the line). Finally, Section A4 will be devoted to some properties of the ⊞-stable distributions with respect to free multiplicative convolution.

## A1. A formula for the density of the ⊞-stable distributions

Case (1) above, namely $\phi(z) = z^{-1}$, corresponds to the semicircle distribution. We start with case (2), so we consider the function

$$\phi_{\alpha,\rho}(z) = e^{i(\alpha-2)\rho\pi} z^{-\alpha+1}$$

with $1 < \alpha < 2$, $0 \le \rho \le 1$.

A1.1. LEMMA. *Let $1 < \alpha < 2$ and $0 \le \rho \le 1$ and define the following region*:

$$\Omega_{\alpha,\rho} = \left\{ re^{-i\theta} \in \mathbf{C}^- \,:\, 0 < \theta < \pi, \quad 0 < r^\alpha < \frac{\sin\theta}{\sin[(2-\alpha)\rho\pi + (\alpha-1)\theta]} \right\}.$$

*Then the map $G_{\nu_{\alpha,\rho}}$ is a one-to-one conformal transformation from $\mathbf{C}^+$ onto $\Omega_{\alpha,\rho}$.*

*Proof.* By a straightforward computation, the region $\Omega_{\alpha,\rho}$ is exactly the set of $z \in \mathbf{C}^-$ such that $\Im(1/z + \phi_{\alpha,\rho}(1/z)) < 0$ so that the lemma follows from [4, Prop. 5.12 (ii)]. □

The region $\Omega_{\alpha,\rho}$ is a Jordan domain and hence, by Carathéodory's theorem (see e.g. [15]), we know that the map $G_{\nu_{\alpha,\rho}}$ extends continuously and gives a homeomorphism of $\mathbf{C}^+ \cup \mathbf{R} \cup \{\infty\}$ with $\overline{\Omega_{\alpha,\rho}}$. We denote this extension again by $G_{\nu_{\alpha,\rho}}$. By the inversion formula for the Cauchy transform, the measure $\nu_{\alpha,\rho}$ has a density with respect to Lebesgue measure which is given by the formula

$$\psi_{\alpha,\rho}(x) = -\frac{1}{\pi}\Im(G_{\nu_{\alpha,\rho}}(x)).$$

For $\rho \in (0,1)$, the boundary of the region $\Omega_{\alpha,\rho}$ consists in the union of the curve given in polar coordinates $z = r e^{-i\theta}$, by the equation

$$r = \left(\frac{\sin\theta}{\sin[(2-\alpha)\rho\pi + (\alpha-1)\theta]}\right)^{\frac{1}{\alpha}}, \quad 0 < \theta < \pi,$$

and the point $0 = G_{\nu_{\alpha,\rho}}(\infty)$, so that for any $x \in \mathbf{R}$, one has $-\Im(G_{\nu_{\alpha,\rho}}(x))/\pi > 0$. Since the function $z \mapsto 1/z + \phi_{\alpha,\rho}(1/z)$ is analytic on $\mathbf{C}^+$ with non-vanishing derivative on $\partial\Omega_{\alpha,\rho}$, it follows that $G_{\nu_{\alpha,\rho}}$ is analytic on $\mathbf{R}$. Let us give an expression for this density. We have seen that for all $x \in \mathbf{R}$, there is a unique $z \in \partial\Omega_{\alpha,\rho}$ such that $1/z + \phi_{\alpha,\rho}(1/z) = x$. Let us call $z(x) = G_{\nu_{\alpha,\rho}}(x) = r(x)e^{-i\theta(x)}$ this point, then one has

$$r(x) = \left(\frac{\sin\theta(x)}{\sin[(2-\alpha)\rho\pi + (\alpha-1)\theta(x)]}\right)^{\frac{1}{\alpha}},$$

and

$$x = \Re\left[\frac{1}{z(x)} + \phi_{\alpha,\rho}\left(\frac{1}{z(x)}\right)\right]$$
$$= (\sin\theta(x))^{-\frac{1}{\alpha}} (\sin[(2-\alpha)\rho\pi + (\alpha-1)\theta(x)])^{\frac{1}{\alpha}-1} \sin[(2-\alpha)\rho\pi + \alpha\theta(x)].$$

The density of the measure $\nu_{\alpha,\rho}$ at the point $x$ is equal to $-\Im z(x)/2\pi$ which is equal to

$$\psi_{\alpha,\rho}(x) = \frac{1}{\pi}(\sin\theta(x))^{1+\frac{1}{\alpha}}(\sin[(2-\alpha)\rho\pi + (\alpha-1)\theta(x)])^{-\frac{1}{\alpha}}.$$

When $\rho = 0$, the boundary of the region $\Omega_{\alpha,\rho}$ is the union of the curve, in polar coordinates $z = re^{-i\theta}$,

$$r = \left(\frac{\sin\theta}{\sin[(2-\alpha)\rho\pi + (\alpha-1)\theta]}\right)^{\frac{1}{\alpha}}, \quad 0 < \theta < \pi,$$

and of the interval $[0, (\alpha-1)^{-1/\alpha}] \subset \mathbf{R}$. The image by $z \mapsto 1/z + \phi_{\alpha,\rho}(1/z)$ of $\partial\Omega_{\alpha,\rho} \cap \mathbf{C}^-$ is the interval $(-\infty, \alpha(\alpha-1)^{1/\alpha-1})$ which is thus the support of the measure $\nu_{\alpha,0}$. On this support, the density is again analytic and given by the formula

$$\psi_{\alpha,\rho}(x) = \frac{1}{\pi}(\sin\theta(x))^{1+\frac{1}{\alpha}}(\sin[(2-\alpha)\rho\pi + (\alpha-1)\theta(x)])^{-\frac{1}{\alpha}}.$$

As remarked earlier, the measure $\nu_{\alpha,1}$ is the image of $\nu_{\alpha,0}$ by the map $t \mapsto -t$ on $\mathbf{R}$. We can summarize the above discussion as follows.

A1.2. PROPOSITION. *For $1 < \alpha < 2$ and $0 < \rho < 1$, the measure $\nu_{\alpha,\rho}$ has a positive analytic density on $\mathbf{R}$ given by*

$$\psi_{\alpha,\rho}(x) = \frac{1}{\pi}(\sin\theta)^{1+\frac{1}{\alpha}}(\sin[(2-\alpha)\rho\pi + (\alpha-1)\theta])^{-\frac{1}{\alpha}},$$

where $\theta \in (0, \pi)$ is the only solution of the equation

$$x = (\sin \theta)^{-\frac{1}{\alpha}} (\sin[(2-\alpha)\rho\pi + (\alpha-1)\theta])^{\frac{1}{\alpha}-1} \sin[(2-\alpha)\rho\pi + \alpha\theta].$$

The support of the measure $\nu_{\alpha,0}$ is the interval $(-\infty, \alpha(\alpha-1)^{1/\alpha-1}]$ and $\nu_{\alpha,0}$ has an analytic density on this support, given by the formula

$$\psi_{\alpha,0}(x) = \frac{1}{\pi} (\sin \theta)^{1+\frac{1}{\alpha}} (\sin(\alpha-1)\theta)^{-\frac{1}{\alpha}},$$

where $\theta \in (0, \pi)$ is the only solution of the equation

$$x = (\sin \theta)^{-\frac{1}{\alpha}} (\sin(\alpha-1)\theta)^{\frac{1}{\alpha}-1} \sin \alpha\theta.$$

Finally, $\nu_{\alpha,1}$ is the image of $\nu_{\alpha,0}$ by the map $t \mapsto -t$.

Similar considerations allow one to treat the cases (3)(ii) and (4). We will state the result below and leave the proof to the reader.

A1.3. PROPOSITION. *For $0 < \rho < 1$, the measure $\nu_{1,\rho}$ has a positive analytic density on $\mathbf{R}$ given by*

$$\psi_{1,\rho}(x) = \frac{1}{\pi} \frac{\sin^2 \theta}{2\rho + \frac{2}{\pi}(1-2\rho)\theta},$$

*where $\theta \in (0, \pi)$ is the only solution of the equation*

$$x = \left(2\rho + \frac{2}{\pi}(1-2\rho)\theta\right) \cot \theta + \frac{2}{\pi}(2\rho-1) \log \frac{\sin \theta}{2\rho + \frac{2}{\pi}(1-2\rho)\theta}.$$

*The support of the measure $\nu_{1,0}$ is the interval $(-\infty, 2/\pi(1+\log \pi/2)]$ and $\nu_{1,0}$ has an analytic density on this support, given by*

$$\psi_{1,0}(x) = \frac{\sin^2 \theta}{2\theta},$$

*where $\theta \in (0, \pi)$ is the only solution of the equation*

$$x = \frac{2}{\pi}\theta \cot \theta + \frac{2}{\pi} \log \frac{\pi \sin \theta}{2\theta}.$$

*The measure $\nu_{1,1}$ is the image of $\nu_{1,0}$ by the map $t \mapsto -t$.*

A1.4. PROPOSITION. *For $0 < \alpha < 1$ and $0 < \rho < 1$, the measure $\nu_{\alpha,\rho}$ has a positive analytic density on $\mathbf{R}$ given by*

$$\psi_{\alpha,\rho}(x) = \frac{1}{\pi} (\sin \theta)^{1+\frac{1}{\alpha}} (\sin[\alpha\rho\pi + (1-\alpha)\theta])^{-\frac{1}{\alpha}},$$

*where $\theta \in (0, \pi)$ is the only solution of the equation*

$$x = (\sin \theta)^{-\frac{1}{\alpha}} (\sin[\alpha\rho\pi + (1-\alpha)\theta])^{\frac{1}{\alpha}-1} \sin[\alpha\rho\pi - \alpha\theta].$$

*The support of the measure $\nu_{\alpha,0}$ is the interval $(-\infty, -\alpha(1-\alpha)^{1/\alpha-1}]$ and $\nu_{\alpha,0}$ has an analytic density on this support given by the formula*

$$\psi_{\alpha,0}(x) = \frac{1}{\pi}(\sin\theta)^{1+\frac{1}{\alpha}}(\sin(1-\alpha)\theta)^{-\frac{1}{\alpha}},$$

*where $\theta \in (0,\pi)$ is the only solution of the equation*

$$x = -(\sin\theta)^{-\frac{1}{\alpha}}(\sin(1-\alpha)\theta)^{\frac{1}{\alpha}-1}\sin\alpha\theta.$$

*The measure $\nu_{\alpha,1}$ is the image of $\nu_{\alpha,0}$ by the map $t \mapsto -t$.*

*Remark.* In the case $\alpha = 1/2$, $\rho = 1$, the density has a closed form, as it is easy to see. Indeed, the support of the measure is $[1/4, +\infty)$, and the density $\psi_{1/2,1}(x) = \sqrt{4x-1}/(2\pi x^2)$.

## A2. Asymptotic behavior and unimodality of the density

The explicit formulas of Section A1 allow us now to derive some properties of ⊞-stable distributions which are reminiscent of properties of ∗-stable distributions. Note that the proof of these properties is much easier than in the classical case.

A2.1. PROPOSITION. *For $1 < \alpha < 2$ and $0 \leq \rho \leq 1$,*

$$\psi_{\alpha,\rho}(x) \sim \frac{1}{\pi}x^{-\alpha-1}\sin(2-\alpha)\rho\pi \quad \text{as } x \to +\infty,$$

$$\psi_{\alpha,\rho}(x) \sim \frac{1}{\pi}(-x)^{-\alpha-1}\sin(2-\alpha)(1-\rho)\pi \quad \text{as } x \to -\infty.$$

*For $0 \leq \rho \leq 1$,*

$$\psi_{1,\rho}(x) \sim \frac{2\rho}{\pi x^2} \quad \text{as } x \to +\infty,$$

$$\psi_{1,\rho}(x) \sim \frac{2(1-\rho)}{\pi x^2} \quad \text{as } x \to -\infty.$$

*For $0 < \alpha < 1$ and $0 \leq \rho \leq 1$,*

$$\psi_{\alpha,\rho}(x) \sim \frac{1}{\pi}x^{-\alpha-1}\sin\alpha\rho\pi \quad \text{as } x \to +\infty,$$

$$\psi_{\alpha,\rho}(x) \sim \frac{1}{\pi}(-x)^{-\alpha-1}\sin\alpha(1-\rho)\pi \quad \text{as } x \to -\infty.$$

*In all these formulas, equivalence to zero means that the function is actually identically zero in a neighborhood of $\pm\infty$.*

*Proof.* Let us check the case where $1 < \alpha < 2$, $0 < \rho < 1$, and $x \to +\infty$. One has then $\theta \to 0$ and $x \sim (\sin\theta)^{-1/\alpha}(\sin(2-\alpha)\rho\pi)^{1/\alpha}$; hence $\sin\theta \sim x^{-\alpha}\sin(2-\alpha)\rho\pi$. Inserting inside the formula for $\psi_{\alpha,\rho}(x)$ gives the result. □

The other cases are similar and left to the reader.

A2.2. PROPOSITION. *The stable distributions are unimodal.*

*Proof.* Recall that a probability distribution on **R**, with a density with respect to Lebesgue measure, is called unimodal if this density has a unique local maximum. Again we will check the result in the case $1 < \alpha < 2$, $0 < \rho < 1$ and leave the other cases to the reader. The function

$$\theta \mapsto (\sin\theta)^{-\frac{1}{\alpha}}(\sin[(2-\alpha)\rho\pi + (\alpha-1)\theta])^{\frac{1}{\alpha}-1}\sin[(2-\alpha)\rho\pi + \alpha\theta]$$

is an increasing diffeomorphism of $(0,\pi)$ onto **R**, so in view of the formula for $\psi_{\alpha,\rho}$ in Proposition A1.2, it is enough to prove that the function

$$\theta \mapsto \frac{1}{\pi}(\sin\theta)^{1+\frac{1}{\alpha}}(\sin[(2-\alpha)\rho\pi + (\alpha-1)\theta])^{-\frac{1}{\alpha}}$$

is unimodal, i.e., has a unique maximum on $(0,\pi)$. Taking the logarithmic derivative and equating to zero we obtain the equation

$$\left(1 + \frac{1}{\alpha}\right)\cot\theta = \left(1 - \frac{1}{\alpha}\right)\cot[(2-\alpha)\rho\pi + (\alpha-1)\theta].$$

Taking the inverse of the cot function, with values in $(0,\pi)$, we obtain the equation

$$\cot^{-1}\left(\frac{\alpha-1}{\alpha+1}\cot[(2-\alpha)\rho\pi + (\alpha-1)\theta]\right) - \theta = 0.$$

A straightforward computation shows that the derivative of this function is strictly increasing for $\theta \in (0,\pi)$, so that the function is strictly convex. Since it is positive at 0 and negative at $\pi$, it must have only one zero in the interval $(0,\pi)$, hence the function $(\sin\theta)^{1+1/\alpha}(\sin[(2-\alpha)\rho\pi + (\alpha-1)\theta])^{-1/\alpha}$ has a unique maximum on $(0,\pi)$. □

### A3. Duality laws

The duality laws for $*$-stable distributions relate the strictly $*$-stable distributions with stability index $\alpha \geq 1$ to the ones with stability index $1/\alpha$. They are discussed in [23, §2.3]. It turns out that the $\boxplus$-stable distributions obey a law of the same form as the $*$-stable distributions, as we show in the next proposition. We will deal only with the case $1 < \alpha \leq 2$ since, as remarked earlier, the strictly $\boxplus$-stable distributions of stability index 1 are the same as the strictly $*$-stable distributions, and in any case the duality law is easy to check in this case.

A3.1. PROPOSITION. *Let $1 < \alpha < 2$ and $\rho \in [0,1]$, and set $\alpha' = 1/\alpha$ and $\rho' = 1 - (2-\alpha)\rho$. Then for all $x > 0$,*

$$\psi_{\alpha,\rho}(x) = x^{-1-\alpha}\psi_{\alpha',\rho'}(x^{-\alpha}).$$

*Proof.* Assume that $\alpha > 1$ and $1 > \rho > 0$. For $x > 0$ consider the unique $\theta \in (0, \pi)$ such that

$$x = (\sin \theta)^{-\frac{1}{\alpha}} (\sin[(2 - \alpha)\rho \pi + (\alpha - 1)\theta])^{\frac{1}{\alpha} - 1} \sin[(2 - \alpha)\rho \pi + \alpha \theta].$$

One has $\alpha \theta + (2 - \alpha)\rho \pi < \pi$, and

$$\psi_{\alpha,\rho}(x) = \frac{1}{\pi}(\sin \theta)^{1+\frac{1}{\alpha}}(\sin[(2 - \alpha)\rho \pi + (\alpha - 1)\theta])^{-\frac{1}{\alpha}}.$$

Note that

$$x^{-\alpha} = \sin \theta \, (\sin[(2 - \alpha)\rho \pi + (\alpha - 1)\theta])^{\frac{1}{\alpha'} - 1}(\sin[(2 - \alpha)\rho \pi + \alpha \theta])^{-\frac{1}{\alpha'}}.$$

With $\theta' = \rho' \pi - \alpha \theta \in (0, \pi)$, we see that

$$\sin \theta = \sin[\alpha' \rho' \pi - \alpha' \theta'],$$
$$\sin[(2 - \alpha)\rho\pi + (\alpha - 1)\theta] = \sin[\alpha' \rho' \pi + (1 - \alpha')\theta'],$$
$$\sin[(2 - \alpha)\rho\pi + \alpha \theta] = \sin \theta'.$$

It follows that

$$x^{-\alpha} = (\sin \theta')^{-\frac{1}{\alpha'}}(\sin[\alpha' \rho' \pi + (1 - \alpha')\theta'])^{\frac{1}{\alpha'} - 1} \sin[\alpha' \rho' \pi - \alpha' \theta'];$$

hence Proposition A1.4 implies that

$$\psi_{\alpha',\rho'}(x^{-\alpha}) = \frac{1}{\pi}(\sin \theta')^{1+\frac{1}{\alpha'}}(\sin[\alpha' \rho' \pi + (1 - \alpha')\theta'])^{-\frac{1}{\alpha'}}.$$

The proposition now follows by a straightforward computation. We leave the case where $\rho = 0$ or $1$ to the reader. □

As in the case of $*$-stable distributions, the duality law can be expressed as an identity in distribution between two random variables, namely

$$Z(\alpha, \rho) \sim Z^{-\frac{1}{\alpha}}\left(\frac{1}{\alpha}, \rho'\right),$$

where $Z(\alpha, \rho)$ denotes a random variable distributed as $\nu_{\alpha,\rho}$, conditioned to be positive; see [23, §3.2, especially Thm. 3.2.5].

Let us also remark that the duality law extends to the case $\alpha = 2$ in the following form. The measure $\nu_{1/2,1}$ is the image of the semicircle distribution by the map $t \mapsto 1/t^2$. This is easily checked using the explicit formula in the remark following Proposition A1.4.

### A4. Multiplicativity properties of ⊞-stable distributions

We will see that the ⊞-stable distributions supported on $(0, +\infty)$ have a nice behavior with respect to the free multiplicative convolution of measures. First recall the definition of the $\Sigma$-transform of a probability distribution $\mu$ on $(0, +\infty)$. Define

$$\Psi_\mu(z) = \int_0^\infty \frac{zt}{1-zt}\, d\mu(t) = \frac{1}{z} G_\mu\left(\frac{1}{z}\right) - 1.$$

Then the function $\Psi_\mu/(1 + \Psi_\mu)$ is univalent on $i\mathbf{C}^+$, with values in $i\mathbf{C}^+$, and we let $\tilde{\chi}_\mu$ denote the inverse function. The $\Sigma$-transform of $\mu$ is the function $\Sigma_\mu(z) = \tilde{\chi}_\mu(z)/z$. It is defined in a neighborhood of $(-\infty, 0)$ in $i\mathbf{C}^+$, and takes positive values on $(-\infty, 0)$. If $\mu$ and $\nu$ are probability measures on $(0, +\infty)$, then there exists a unique probability measure on $(0, +\infty)$, denoted $\mu \boxtimes \nu$, such that

$$\Sigma_{\mu \boxtimes \nu} = \Sigma_\mu \Sigma_\nu$$

in a neighborhood of some interval $(-\varepsilon, 0)$ in $\mathbf{C}$, where these functions are defined. The measure $\mu \boxtimes \nu$ is called the *free multiplicative convolution* of the measures $\mu$ and $\nu$. This binary operation was introduced by Voiculescu in [18] (see also [3] and [4]).

According to Propositions A1.2, A1.3, and A1.4, the only ⊞-stable distributions which have their support in $(0, +\infty)$ are the measures $\nu_{\alpha,1}$ with $0 < \alpha < 1$. For $0 < \alpha < 1$, let $\nu_\alpha = \nu_{\alpha,1}$ and $\sigma_\alpha$ be the positive one-sided $*$-stable distribution of stability index $\alpha$, so that

$$\int_0^\infty e^{-xt}\, d\sigma_\alpha(t) = e^{-x^\alpha}, \quad x \geq 0.$$

For a probability measure $\mu$ on $(0, +\infty)$, we will denote $\check{\mu}$ the image of $\mu$ by the map $t \mapsto 1/t$ on $(0, \infty)$.

A4.1. LEMMA. *For all $\alpha \in (0, 1)$,*

$$\Sigma_{\nu_\alpha}(z) = \left(\frac{-z}{1-z}\right)^{\frac{1}{\alpha}-1} \quad \text{and} \quad \Sigma_{\check{\nu}_\alpha}(z) = (1-z)^{\frac{1}{\alpha}-1}.$$

*Proof.* As noted above, $\Psi_\mu(z) = G_\mu(1/z)/z - 1$ for any probability measure $\mu$ on $(0, +\infty)$ and $z \in \mathbf{C} \setminus \mathbf{R}^+$. Since

$$\frac{1}{G_{\nu_\alpha}(z)} - e^{i\alpha\pi} G_{\nu_\alpha}(z)^{\alpha-1} = \frac{1}{rG_{\nu_\alpha}(z)}(1 - (-G_{\nu_\alpha}(z))^\alpha) = z,$$

for $z \in (-\infty, 0)$, we see that

$$\frac{1}{z(1 + \Psi_{\nu_\alpha}(z))}(1 - (-z(1 + \Psi_{\nu_\alpha}(z)))^\alpha) = \frac{1}{z};$$

hence

(A4.2) $\qquad -\Psi_{\nu_\alpha}(z) = (-z(1+\Psi_{\nu_\alpha}(z)))^\alpha = (-z)^\alpha(1+\Psi_{\nu_\alpha}(z))^\alpha$

for $z \in (-\infty, 0)$, since then $1 + \Psi_{\nu_\alpha}(z) \in (0, +\infty)$. Solving this equation with respect to $\zeta = \Psi_{\nu_\alpha}(z)/(1+\Psi_{\nu_\alpha}(z)) \in (-\infty, 0)$, we get

$$-z = (-\zeta)^{1/\alpha}(1-\zeta)^{1-1/\alpha};$$

hence

$$\tilde{\chi}_{\nu_\alpha}(\zeta) = -(-\zeta)^{\frac{1}{\alpha}}(1-\zeta)^{1-\frac{1}{\alpha}},$$

and

$$\Sigma_{\nu_\alpha}(\zeta) = \left(\frac{-\zeta}{1-\zeta}\right)^{\frac{1}{\alpha}-1}, \quad \zeta \in (-\infty, 0).$$

The claim follows by analytic continuation. $\qquad \square$

For any probability measure $\mu$ on $(0, +\infty)$, $\Psi_\mu(z) + \Psi_{\check{\mu}}(1/z) + 1 = 0$ where $\check{\mu}$ is the image of $\mu$ by $t \mapsto 1/t$. From (A4.2), it follows that

$$1 + \Psi_{\check{\nu}_\alpha}(z) = \left(\frac{1}{z}\Psi_{\check{\nu}_\alpha}(z)\right)^{\alpha-1}$$

for all $z \in (-\infty, 0)$; hence solving with respect to $\zeta = \Psi_{\check{\nu}_\alpha}(z)/(1+\Psi_{\check{\nu}_\alpha}(z))$ we get

$$-\tilde{\chi}_{\check{\nu}_\alpha}(\zeta) = (-\zeta)(1-\zeta)^{\frac{1}{\alpha}-1}$$

and

$$\Sigma_{\check{\nu}_\alpha}(\zeta) = (1-\zeta)^{\frac{1}{\alpha}-1} \qquad \text{for} \quad \zeta \in (-\infty, 0).$$

Again this gives the result by analytic continuation.

A4.3. PROPOSITION. *For all $s, t > 0$,*

$$\nu_{\frac{1}{1+s}} \boxtimes \nu_{\frac{1}{1+t}} = \nu_{\frac{1}{1+s+t}}.$$

*Proof.* This is a straightforward consequence of Lemma A4.1, and the multiplicativity of the $\Sigma$ transform with respect to free multiplicative convolution. $\qquad \square$

We turn now to a remarkable relation between the $\boxplus$-stable distributions and the $*$-stable distributions.

A4.4. PROPOSITION. *For every $\alpha \in (0, 1)$,*

$$\nu_\alpha \boxtimes \check{\nu}_\alpha = \sigma_\alpha \circledast \check{\sigma}_\alpha,$$

*where $\circledast$ denotes the convolution of measures on the multiplicative group $(0, +\infty)$.*

*Proof.* From Lemma A4.1, $\Sigma_{\nu_\alpha \boxtimes \check{\nu}_\alpha}(z) = (-z)^{1/\alpha - 1}$. This $\Sigma$-transform can be inverted (see [5, §5.4]), and the measure $\nu_\alpha \boxtimes \check{\nu}_\alpha$ turns out to have density

$$\frac{\frac{1}{\pi}\sin(\alpha\pi)y^{\alpha-1}}{y^{2\alpha} + 2\cos(\alpha\pi)y^\alpha + 1}, \quad y > 0,$$

with respect to Lebesgue measure. As remarked in [5, §5.4], this measure is also equal to $\sigma_\alpha \circledast \check{\sigma}_\alpha$. □

INDIANA UNIVERSITY, BLOOMINGTON, IN
*E-mail address*: bercovic@indiana.edu

DIPARTIMENTO DI MATEMATICA, UNIVERSITÀ DI BRESCIA, BRESCIA, ITALY
*E-mail address*: pata@bsing.unibs.it

CNRS, UNIV. PIERRE ET MARIE CURIE, PARIS, FRANCE AND
DMI, ÉCOLE NORMALE SUPÉRIEURE, PARIS, FRANCE
*E-mail address*: pbi@ccr.jussieu.fr

## REFERENCES

[1] H. BERCOVICI and V. PATA, Classical versus free domains of attraction, Math. Res. Lett. **2** (1995), 791–795.
[2] ______, The law of large numbers for free identically distributed random variables, Ann. Probab. **24** (1996), 453–465.
[3] H. BERCOVICI and D. VOICULESCU, *Lévy-Hinčin type theorems for multiplicative and additive free convolution*, Pacific J. Math. **153** (1992), 217–248.
[4] ______, Free convolution of measures with unbounded support, Indiana Univ. Math. J. **42** (1993), 733–773.
[5] P. BIANE, Processes with free increments, Math. Z. **227** (1998), 143–174.
[6] N. H. BINGHAM, C. M. GOLDIE, and J. L. TEUGELS, *Regular Variation*, John Wiley & Sons, New York, 1971.
[7] W. FELLER, *An Introduction to Probability Theory and its Applications*, Cambridge University Press, Cambridge, 1987.
[8] B. V. GNEDENKO and A. N. KOLMOGOROV, *Limit Distributions for Sums of Independent Random Variables*, Addison-Wesley Publ. Co., Cambridge, Mass., 1954.
[9] P. LÉVY, *Théorie de L'addition des Variables Aléatoires*, Gauthier-Villars, Paris, 1937.
[10] J. M. LINDSAY and V. PATA, Some weak laws of large numbers in noncommutative probability, Math. Z. **226** (1997), 533–543.
[11] H. MAASSEN, Addition of freely independent random variables, J. Funct. Anal. **106** (1992), 409–438.
[12] V. PATA, Lévy type characterization of stable laws for free random variables, Trans. A.M.S. **347** (1995), 2457–2472.
[13] ______, Domains of partial attraction in noncommutative probability, Pacific J. Math. **176** (1996), 235–248.
[14] ______, The central limit theorem for free additive convolution, J. Funct. Anal. **140** (1996), 359–380.
[15] W. RUDIN, *Real and Complex Analysis*, McGraw-Hill Book Co., New York, 1966.
[16] R. SPEICHER and R. WOROUDI, Boolean convolution, Fields Inst. Commun. **12** (1997), 267–279.

[17] D. Voiculescu, Symmetries of some reduced free product C*-algebras, *Operator Algebras and their Connections with Topology and Ergodic Theory*, Lecture Notes in Math. **1132**, Springer-Verlag, New York, 556–588.
[18] ———, Addition of certain noncommuting random variables, J. Funct. Anal. **66** (1986), 323–346.
[19] ———, Multiplication of certain noncommuting random variables, J. Operator Theory **18** (1987), 223–235.
[20] ———, Limit laws for random matrices and free products, Invent. Math. **104** (1991), 201–220.
[21] D. Voiculescu, K. Dykema, and A. Nica, Free random variables, CRM Monograph Series, No. 1, A.M.S., Providence, RI, 1992.
[22] E. P. Wigner, On the distribution of the roots of certain symmetric matrices, Ann. of Math. **67** (1958), 325–327.
[23] V. M. Zolotarev, *One-dimensional Stable Distributions*, Transl. of Math. Monographs **65**, A.M.S., Providence, RI, 1986.